\documentclass{article}
\usepackage[T2A]{fontenc}
\usepackage[cp1251]{inputenc}
\usepackage[english,russian]{babel}
\usepackage[tbtags]{amsmath}
\usepackage{amsfonts,amssymb,mathrsfs,amscd}
\numberwithin{equation}{section}
\newtheorem{remark}{Remark}[section]
\newtheorem{lemma}{Lemma}[section]
\newtheorem{theorem}{Theorem}[section]

\newtheorem{conclusion}{Conclusion}

\newcommand{\res}{\mathop{\rm res}}

\newcommand{\col}{\mathop{\rm col}}

\begin{document}
\begin{Large}
\thispagestyle{empty}
\begin{center}
{\bf Inverse scattering problem for third order differential operators on the whole axis\\
\vspace{5mm}
V. A. Zolotarev}\\

B. Verkin Institute for Low Temperature Physics and Engineering
of the National Academy of Sciences of Ukraine\\
47 Nauky Ave., Kharkiv, 61103, Ukraine

Department of Higher Mathematics and Informatics, V. N. Karazin Kharkov National University \\
4 Svobody Sq, Kharkov, 61077,  Ukraine

\end{center}
\vspace{5mm}

{\small {\bf Abstract.} Inverse scattering problem for an operator, which is a sum of the operator of the third derivative and of an operator of multiplication by a real function, is solved. The main closed system of equations of inverse problem is obtained. This system contains scattering coefficients and bound state elements as independent parameters. Form of the simplest reflectionless potential (analogous to soliton) is found.}
\vspace{5mm}


\vspace{5mm}

\begin{center}
{\bf Introduction}
\end{center}
\vspace{5mm}

Scattering theory is one of the interesting, intriguing parts of mathematical physics and functional analysis. A key role in this theory is played by the method of inverse scattering problem which became one of the main instrument for solving linear equations in partial derivatives. In one-dimensional case, for operators of second order (Sr\"{o}dinger operator, in particular), method of the inverse problem was developed by V. A. Marchenko, I. M. Gelfand -- B. M. Levitan, M. G. Kreyn, L. D. Faddeev \cite{1} -- \cite{4}. Numerous applications of the method of inverse scattering problem, as well as its history, are given in the book by K. Chadan, P. C. Sabatier.

Generalization of the method of inverse problem for differential operators of order $n$ ($n\geq3$) didn't give the wanted results, apart from particular cases, which led to the narrowing of the studied class of operators or to additional conditions of similar nature \cite{6} -- \cite{9}. Probably, this is connected with the absence of transformation operators for $n>2$. Integration of non-linear equations, i.e., construction of L - A-pairs for these equations (e.g., for the Degasperis -- Procesi equation \cite{10} -- \cite{14}) leads to a third order operator $L$ (cubic string). Therefore, there exists an urgent need in creation of a method of inverse scattering problem for differential operators of $n$th order (and for $n=3$ in particular).

In works \cite{15} -- \cite{18}, not only inverse problems for third order differential operators (spectral and scattering) are solved, but also the needed for study of such problems analytical apparatus is created. This apparatus is easily generalized for operators of arbitrary order.

This work is dedicated to the scattering problem for the operator
\begin{equation}
(L_qy)(x)=iy'''(x)+q(x)y(x)\label{eq0.1}
\end{equation}
on the whole axis, $x\in\mathbb{R}$, where potential $q(x)$ is real and satisfies the condition
\begin{equation}
\int\limits_{-\infty}^\infty|q(x)|^2e^{2a|x|}dx<\infty\label{eq0.2}
\end{equation}
($a>0$ is a fixed number). The paper consists of three sections, the first of which deals with Jost solutions. In this case, we have three Jost solutions that are transformed into each other after the turn at an angle $2\pi/3$ (an analogue of $\lambda\rightarrow-\lambda$ for Sturm -- Liouville). Each of these solutions is holomorphic relative to $\lambda$ in one of the three sectors of half-plane $\mathbb{C}$. Asymptotic properties of Jost solutions in these sectors are described.

The second section is dedicated to the scattering problem. As usual \cite{3}, an important role here is played by the matrix of transition from one system of fundamental solutions $\{u_k(\lambda,x)\}_0^2$ at ``$-\infty$'' to another system $\{v_k(\lambda,x)\}_0^2$ at ``$+\infty$''. It is proved that this matrix is $J$-unitary and unimodular. For the Sr\"{o}dinger operator, transition from one system of solutions to another is conducted by involution, viz., ``taking the complex conjugate'', whereas in our case, transition to the dual system takes place as Wronskian is computed (Lie bracket). Presence of three waves at ``$-\infty$'' and of three waves at ``$+\infty$'' leads to the scattering problem in which there are two scattering coefficients and one transmission coefficient. It is proved that the numbers corresponding to bound states are situated on rays in a sector, some of which correspond to the positive eigenvalues and others, to the negative eigenvalues. Conditions for ``$a$'' from $[0,2)$ when the number of bound states is finite are found. Using Riemann boundary value problems, problems on jumps on the borders of sectors of Jost solutions holomorphicity are derived. A closed system of singular equations (an analogue of Marchenko equation) for the solution to the inverse problem is found. Scattering coefficients and bound state parameters enter this system as independent parameters.

The last, third, section solves the inverse scattering problem. The explicit form of potential is found in terms of solutions to the main equation system. Form of the simplest reflectionless potential (a soliton analogue) is found.

\section{Jost solutions}\label{s1}

{\bf 1.1} Equation
\begin{equation}
iy'''(x)=\lambda^3y(x)\quad(x\in\mathbb{R},\lambda\in\mathbb{C})\label{eq1.1}
\end{equation}
has three linear independent solutions $\{e^{i\lambda\zeta_kx}\}_0^2$ \cite{15,16,9}, where $\zeta_k$ are roots of the equation $\zeta^3=1$, and
\begin{equation}
\zeta_0=1;\quad\zeta_1=\frac12(-1+i\sqrt3);\quad\zeta_2=\frac12(-1-i\sqrt3),\label{eq1.2}
\end{equation}
and any solution to \eqref{eq1.1} is their linear combination.  Instead of $\{e^{i\lambda\zeta_kx}\}_0^2$, it is convenient to take the system of fundamental solutions $\{s_p(i\lambda\zeta_k)\}_0^2$ \cite{15} -- \cite{18},
\begin{equation}
s_p(z)=\frac13\sum\limits_{k=0}^2\frac1{\zeta_k^p}e^{\zeta_kz}\quad(0\leq p\leq2;z\in\mathbb{C}).\label{eq1.3}
\end{equation}
The functions $\{s_p(i\lambda x)\}_0^2$ play an important role in the study of third order differential operators \cite{15} -- \cite{18}, and their main properties are as follows.

\begin{lemma}[\cite{15} -- \cite{18}]\label{l1.1}
Entire functions of exponential type $\{s_p(z)\}_0^2$ \eqref{eq1.3} have the following properties:

(i) $s'_p(z)=s_{p'}(z)$ ($p'=(p+2)\mod3$; $y'(z)=dy/dz$);

(ii) $\overline{s_p(z)}=s_p(\overline{z})$ ($0\leq p\leq2$);

(iii) {\bf $p$-evenness}, $s_p(z\zeta_1)=\zeta_1^ps_p(z)$ $(0\leq p<2)$;

(iv) Euler formula
$$e^{z\zeta_k}=s_0(z)+\zeta_ks_1(z)+\zeta_k^2s_2(z)\quad(0\leq k\leq2);$$

(v) the functions $\{s_p(z)\}_0^2$ \eqref{eq1.3} are solutions to the equation $y'''(z)=y(z)$ and satisfy the initial data:
$$s_0(0)=1;\quad s'_0(0)=0;\quad s''_0(0)=0;$$
$$s_1(0)=0;\quad s'_1(0)=1;\quad s''_1(0)=0;$$
$$s_2(0)=0;\quad s'_2(0)=0;\quad s''_2(1)=1;$$

(vi) the main identity
$$s_0^3(z)+s_1^3(z)+s_2^3(z)-3s_0(z)s_1(z)s_2(z)=1;$$

(vii) addition formulas
$$s_0(z+w)=s_0(z)s_0(w)+s_1(z)s_2(w)+s_2(z)s_1(w);$$
$$s_1(z+w)=s_0(z)s_1(w)+s_1(z)s_0(w)+s_2(z)s_2(w);$$
$$s_2(z+w)=s_0(z)s_2(w)+s_1(z)s_1(w)+s_2(z)s_0(w);$$

(viii)
$$3s_0(z)s_0(w)=s_0(z+w)+s_0(z+\zeta_1w)+s_0(z+\zeta_2w);$$
$$3s_0(z)s_2(w)=s_0(z+w)+\zeta_2s_0(z+\zeta_1w)+\zeta_2s_0(z+\zeta_2w);$$
$$3s_1(z)s_0(w)=s_1(z+w)+s_1(z+\zeta_1w)+s_1(z+\zeta_2w);$$
$$3s_2(z)s_2(w)=s_1(z+w)+\zeta_1s_1(z+\zeta_1w)+\zeta_2s_1(z+\zeta_2w);$$
$$3s_2(z)s_0(w)=s_2(z+w)+s_2(z+\zeta_1w)+s_2(z+\zeta_2w);$$
$$3s_1(z)s_1(w)=s_2(z+w)+\zeta_2s_2(z+\zeta_1w)+\zeta_1s_2(z+\zeta_2w);$$

(ix)
$$3s_0^2(z)=s_0(2z)+2s_0(-z);$$
$$3s_1^2(z)=s_2(2z)+2s_2(-z);$$
$$3s_2^2(z)=s_1(2z)+2s_1(-z);$$

(x)
$$s_0^2(z)-s_1(z)s_2(z)=s_0(-z);$$
$$s_1^2(z)-s_0(z)s_2(z)=s_2(-z);$$
$$s_2^2(z)-s_0(z)s_1(z)=s_1(-z);$$

(xi)
Taylor formula
$$s_0(z)=1+\frac{z^3}{3!}+\frac{z^6}{6!}+...;$$
$$s_1(z)=z+\frac{z^4}{4!}+\frac{z^7}{7!}+...;$$
$$s_2(z)=\frac{z^2}{2!}+\frac{z^5}{5!}+\frac{z^8}{8!}+....$$
\end{lemma}

Solution to the Cauchy problem
\begin{equation}
iy'''(x)=\lambda^3y(x)+f(x);\quad y(0)=y_0;\quad y'(0)=y_1;\quad y''(0)=y_2;\label{eq1.4}
\end{equation}
for $f=0$, is
\begin{equation}
y_0(\lambda,x)\stackrel{\rm def}{=}y_0s_0(i\lambda x)+y_1\frac{s_1(i\lambda x)}{i\lambda}+y_2\frac{s_2(i\lambda x)}{(i\lambda)^2}.\label{eq1.5}
\end{equation}
Using the variation of constants method, hence we find the solution to the Cauchy problem \eqref{eq1.4}
\begin{equation}
y(\lambda,x)=y_0(\lambda,x)-i\int\limits_0^x\frac{s_2(i\lambda(x-t))}{(i\lambda)^2}f(t)dt\label{eq1.6}
\end{equation}
where $y_0(\lambda,x)$ is given by \eqref{eq1.5} and $s_2(z)$, correspondingly, by \eqref{eq1.3}.

Unit vectors $\{\zeta_k\}_0^2$ \eqref{eq1.2} specify the straight lines
\begin{equation}
L_{\zeta_k}\stackrel{\rm def}{=}\{x\zeta_k:x\in\mathbb{R}\}\quad(0\leq k\leq2),\label{eq1.7}
\end{equation}
and let $l_{\zeta_k}$ be the rays from the origin in the direction of $\zeta_k$ (and $\widehat{l}_{\zeta_k}$ be the rays coming to the origin),
\begin{equation}
l_{\zeta_k}\stackrel{\rm def}{=}\{x\zeta_k:x\in\mathbb{R}_+\};\quad\widehat{l}_{\zeta_k}=L_{\zeta_k}\backslash l_{\zeta_k}\stackrel{\rm def}{=}\{x\zeta_k:x\in\mathbb{R}_-\}\quad(0\leq k\leq2).\label{eq1.8}
\end{equation}
Straight lines $L_{\zeta_k}$ \eqref{eq1.7} divide the plane $\mathbb{C}$ into six sectors:
\begin{equation}
S_p\stackrel{\rm def}{=}\left\{z\in\mathbb{C}:\frac{2\pi}6(p-1)<\arg z<\frac{2\pi}6\pi\right\}\quad(1\leq p\leq6).\label{eq1.9}
\end{equation}
\vspace{5mm}

{\bf 1.2} Consider the self-adjoint operator $L_q$ in space $L^2(\mathbb{R})$,
\begin{equation}
(L_q)(x)\stackrel{\rm def}{=}iy'''(x)+q(x)y(x),\label{eq1.10}
\end{equation}
domain of which is
\begin{equation}
\mathfrak{D}(L_q)\stackrel{\rm def}{=}\{y:y\in W_2^3(\mathbb{R})\},\label{eq1.11}
\end{equation}
and potential is real and
\begin{equation}
\int\limits_{\mathbb R}|q(x)|^2e^{2a|x|}dx<0\quad(a>0).\label{eq1.12}
\end{equation}

\begin{remark}\label{r1.1}
Formula \eqref{eq1.12} implies that $q(x)e^{b|x|}\in L^1(\mathbb{R})\cap L^2(\mathbb{R})$ for all $b$ such that $0\leq b<a$. If $q=0$, then $L_q$ coincides with the unperturbed operator $L_0$ ($=iD^3$, ${\displaystyle D=\frac d{dx}})$.
\end{remark}

Equation
\begin{equation}
iy'''(x)+q(x)y(x)=\lambda^3y(x)\label{eq1.13}
\end{equation}
when $x\rightarrow\pm\infty$, transforms into equation \eqref{eq1.1}, due to \eqref{eq1.12}. Therefore it is natural to define the {\bf Jost solutions} \cite{1} -- \cite{3} which are solutions to equation \eqref{eq1.13} and behave as exponents $\{e^{i\lambda\zeta_kx}\}$ as $x\rightarrow\pm\infty$,
\begin{equation}
\begin{array}{ccc}
(a)\,v_k(\lambda,x)-e^{i\lambda\zeta_kx}\rightarrow0&(x\rightarrow\infty,0\leq k\leq2);\\
(b)\,u_k(\lambda,x)-e^{i\lambda\zeta_kx}\rightarrow0&(x\rightarrow-\infty,0\leq k\leq2).
\end{array}\label{eq1.14}
\end{equation}
Show that the integral equation
\begin{equation}
v_k(\lambda,x)=e^{i\lambda\zeta_kx}-i\int\limits_x^\infty\frac{s_2(\lambda(x-t))}{(i\lambda)^2}q(t)v_k(t)dt\quad(0\leq k\leq2)\label{eq1.15}
\end{equation}
is equivalent to the boundary value problem \eqref{eq1.13}, (a) \eqref{eq1.14}. Determine solvability of equation \eqref{eq1.15}. In $L^2(\mathbb{R})$, define Volterra operator family (depending on $\lambda$),
\begin{equation}
(K_\lambda f)(x)\stackrel{\rm def}{=}\int\limits_x^\infty K_1(\lambda,x,t)q(t)f(t)dt\quad(f\in L^2(\mathbb{R}))\label{eq1.16}
\end{equation}
where
\begin{equation}
K_1(\lambda,x,t)=\frac{s_2(i\lambda(x-t))}{(i\lambda)^2}.\label{eq1.17}
\end{equation}
Equation \eqref{eq1.15} in terms of $K_\lambda$ becomes
$$(I+iK_\lambda)v_k(\lambda,x)=e^{i\lambda\zeta_kx}\quad(0\leq k\leq2),$$
and thus
\begin{equation}
v_k(\lambda)=\sum\limits_0^\infty(-i)^nK_\lambda^ne^{i\lambda\zeta_kx}.\label{eq1.18}
\end{equation}
Operators $K_\lambda^n$ are also Volterra,
$$(K_\lambda^nf)(x)\stackrel{\rm def}{=}\int\limits_x^\infty K_n(\lambda,x,t)q(t)f(t)\quad(f\in L^2(\mathbb{R}),$$
and for the kernels $K_n(\lambda,x,t)$ the following recurrent relations hold:
\begin{equation}
K_{n+1}(\lambda,x,t)=\int\limits_x^tK_n(\lambda,x,s)q(s)K_1(\lambda,s,t)ds\quad(n\in\mathbb{N}).\label{eq1.19}
\end{equation}

\begin{lemma}\label{l1.2}
For the kernels $K_n(\lambda,x,t)$ \eqref{eq1.19}, the following equalities are true:
\begin{equation}
K_n(\lambda\zeta_1,x,t)=K_n(\lambda,x,t);\quad\overline{K_n(\lambda,x,t)}=K_n(\overline{\lambda},t,x);\label{eq1.20}
\end{equation}
and the following estimates hold:
\begin{equation}
|K_n(\lambda,x,t)|\leq\left\{
\begin{array}{lll}
{\displaystyle\frac{d(\lambda(t-x))}{|\lambda|^{2n}}\cdot\frac{\sigma^{n-1}(t)}{(n-1)!}\quad(\lambda\not=0,t\geq x);}\\
{\displaystyle\left(\frac{(t-x)^2}2\right)^n\frac{\sigma^{n-1}(t)}{n^{2n}(n-1)!}\quad(\lambda=0),}
\end{array}\right.\label{eq1.21}
\end{equation}
here
\begin{equation}
d(\lambda)\stackrel{\rm def}{=}e^{|\beta|}\ch\frac{\alpha\sqrt3}2\quad(\lambda=\alpha+i\beta;\alpha,\beta\in\mathbb{R});\quad\sigma(t)=\int\limits_{-\infty}^t|q(s)|ds.\label{eq1.22}
\end{equation}
\end{lemma}

P r o o f. Relations \eqref{eq1.20} follow from the recurrent formulas \eqref{eq1.19} ($q(x)$ is real) and the fact that equalities \eqref{eq1.20} hold for $K_1(\lambda,x,t)$ \eqref{eq1.17}. Since
\begin{equation}
\begin{array}{ccc}
{\displaystyle i\zeta_0\lambda=-\beta+\alpha i;\quad i\zeta_1\lambda=\frac12[(\beta-\alpha\sqrt3)-i(\alpha+\beta\sqrt3)];}\\
{\displaystyle i\zeta_2\lambda=\frac12[(\beta+\alpha\sqrt3)-i(\alpha-\beta\sqrt3)]}
\end{array}\label{eq1.23}
\end{equation}
($\lambda=\alpha+i\beta$; $\alpha$, $\beta\in\mathbb{R}$), then
$$s_k(i\lambda)=\frac13\left\{e^{-\beta+\alpha i}+\frac1{\zeta_1^k}e^{\frac12((\beta-\alpha\sqrt3)+i(\alpha+\beta\sqrt3))}+\frac1{\zeta_2^k}e^{\frac12((\beta+\alpha\sqrt3)-i(\alpha-\beta\sqrt3))}\right\},$$
and thus $|s_k(i\lambda)|\leq d(\lambda)$ ($0\leq k\leq2$) where $d(\lambda)$ is from \eqref{eq1.22}. So, for $K_1(\lambda,x,t)$ \eqref{eq1.17} we have
$$|K_1(\lambda,x,t)|\leq\frac{d(\lambda(t-x))}{|\lambda|^2}\quad(t\geq x),$$
this coincides with \eqref{eq1.21} for $n=1$ ($\lambda\not=0$). Using induction by $n$ and the recurrent relation \eqref{eq1.19}, we obtain
$$|K_{n+1}(\lambda,x,t)|\leq\int\limits_x^t\frac{d(\lambda(s-x))}{|\lambda|^{2n}}\frac{\sigma^{n-1}(s)}{(n-1)!}\cdot|q(s)|\frac{d(\lambda(t-s))}{|\lambda|^2}ds,$$
and since
$$d(\lambda(s-x))d(\lambda(t-s))=e^{|\beta|(t-x)}\ch\frac{\alpha\sqrt3}2(s-x)\ch\frac{\alpha\sqrt3}2(t-s)=\frac12e^{|\beta|(t-x)}$$
$$\times\left[\ch\frac{\alpha\sqrt3}2(t-x)+\ch\frac{\alpha\sqrt3}2(t+x-2s)\right]\leq d(\lambda(t-x)),$$
then
$$|K_{n+1}(\lambda,x,t)|\leq\frac{d(\lambda(t-x))}{|\lambda|^{2(n+1)}}\frac{\sigma^n(t)}{n!}.$$

If $\lambda=0$, then due to \eqref{eq1.17} we have that ${\displaystyle K_1(0,x,t)=\frac{(t-x)^2}2}$, which coincides with \eqref{eq1.21} ($n=1$, $\lambda=0$). Again using induction by $n$ and \eqref{eq1.19}, we obtain
$$|K_{n+1}(0,x,t)|\leq\int\limits_x^t\left(\frac{(t-x)^2}2\right)^n\frac{\sigma^{n-1}(s)}{n^{2n}(n-1)!}|q(s)|\frac{(t-s)^2}2ds.$$
Function $f(s)=(t-s)^2(s-x)^{2n}$ is positive for $s\in(x,t)$ and $f(x)=f(t)=0$. It has maximum on $(x,t)$ at the point $s_0=(x+nt)/n+1$ and $f(s_0)=(t-x)^{2(n+1)}\cdot n^{2n}/(n+1)^{2(n+1)}$, and thus
$$|K_{n+1}(0,x,t)|\leq\left(\frac{(t-x)^2}2\right)^{n+1}\frac1{(n+1)^{2(n+1)}}\frac{\sigma^n(t)}{n!}.\blacksquare$$

Formula \eqref{eq1.18} implies that
\begin{equation}
v_k(\lambda,x)=e^{i\lambda\zeta_kx}+\int\limits_x^\infty N(\lambda,x,t)q(t)e^{i\lambda\zeta_kt}dt\quad(0\leq k\leq2)\label{eq1.24}
\end{equation}
where
\begin{equation}
N(\lambda,x,t)=\sum\limits_1^\infty(-1)^nK_n(\lambda,x,t).\label{eq1.25}
\end{equation}
Series \eqref{eq1.25} is majorized by the convergent (due to \eqref{eq1.21}) series, and
\begin{equation}
|N(\lambda,x,t)|\leq\left\{
\begin{array}{lll}
{\displaystyle\frac{d(\lambda(x-t))}{|\lambda|^2}\cdot\exp\left\{\frac{\sigma(t)}{|\lambda|^2}\right\}\quad(\lambda\not=0);}\\
{\displaystyle\frac{(t-x)^2}2\cdot\sum\limits_0^\infty\left(\frac{(t-x)^2}2\sigma(t)\right)^n\frac1{(n+1)^{(n+1)2}\cdot n!}\quad(\lambda=0).}
\end{array}\right.\label{eq1.26}
\end{equation}
Hence the uniform convergence (in $x$) of integral \eqref{eq1.24} for $|\lambda|<a$ follows, and thus $v_k(\lambda,x)$ is differentiable by $x$. Analogously it is established that $v_k(\lambda,x)$ has third derivative by $x$ and satisfies equation \eqref{eq1.13}.

Proceed to boundary condition (a) \eqref{eq1.14}. Since $x\rightarrow+\infty$, we assume that $x>0$, then \eqref{eq1.24} ($\lambda\not=0$) implies that
$$|v_k(\lambda,x)|\leq e^{|\lambda|x}+\int\limits_x^\infty dt|q(t)|e^{|\lambda|t}\frac{d(\lambda(t-x))}{|\lambda|^2}\exp\left\{\frac{\sigma(t)}{|\lambda|^2}\right\}\leq e^{|\lambda|x}$$
$$+\frac1{|\lambda|^2}\int\limits_x^\infty dt|q(t)|e^{|\lambda|t}e^{\lambda|\lambda|(t-x)}\cdot\exp\left|\frac{\sigma(t)}{|\lambda|^2}\right|\leq e^{|\lambda|x}\left(1+\frac{\exp\left({\displaystyle\frac\sigma{|\lambda|^2}}\right)}{|\lambda|^2}\int\limits_x^\infty dte^{3|\lambda|(t-x)}|q(t)|\right)$$
where $\sigma\stackrel{\rm def}{=}\sigma(\infty)$. Using $e^{3|\lambda|(t-x)}|q(t)|\leq e^{3|\lambda|t}|q(t)|\in L^1(\mathbb{R})$ (when $3|\lambda|<a$, Remark \ref{r1.1}), we obtain that
$$|v_k(\lambda,x)|<e^{|\lambda|\cdot|x|}\left(1+\frac{q_1}{|\lambda|^2}\exp\left\{\frac\sigma{|\lambda|^2}\right\}\right)\quad(x>0,3|\lambda|<a),$$
here $q_1=\|e^{3|\lambda|t}|q(t)|\|_{L^1(\mathbb{R})}$. Hence it follows that
$$\left|\int\limits_x^\infty\frac{s_1(i\lambda(x-t))}{(i\lambda)^3}q(t)v_k(\lambda,t)dt\right|\leq\left(1+\frac{q_1}{|\lambda|^2}\exp\left\{\frac\sigma{|\lambda|^2}\right\}\right)
\int\limits_x^\infty\frac{d(\lambda(t-x))}{|\lambda|^2}|q(t)|e^{|\lambda|t}dt$$
$$\leq\frac1{|\lambda|^2}\left(1+\frac{q_1}{|\lambda|^2}\exp\left\{\frac\sigma{|\lambda|^2}\right\}\right)\int\limits_x^\infty e^{3|\lambda|t}|q(t)|dt\quad(x>0).$$
If $3|\lambda|<a$, then $e^{3|\lambda||t|}|q(t)|\in L^1(\mathbb{R})$, therefore the last integral tends to zero as $x\rightarrow\infty$, this provides meeting of the boundary condition (a) \eqref{eq1.14}, due to \eqref{eq1.15}. Define the circle
\begin{equation}
\mathbb{D}_{a/3}\stackrel{\rm def}{=}\{\lambda\in\mathbb{C}:|\lambda|<a/3\}.\label{eq1.27}
\end{equation}

\begin{theorem}\label{t1.1}
Let a potential $q(x)$ be real and \eqref{eq1.12} hold, then for all $\lambda\in\mathbb{D}_{a/3}$ \eqref{eq1.27} the Jost solutions $\{v_k(\lambda,x)\}_0^2$ of the boundary value problem \eqref{eq1.13}, (a) \eqref{eq1.14} exist and are given by \eqref{eq1.24}. Functions $\{v_k(\lambda,x)\}_0^2$ are holomorphic by $\lambda$ in $\mathbb{D}_{a/3}$ and
\begin{equation}
v_k(\lambda\zeta_1,x)=v_{k'}(\lambda,x)\quad(k'=(k+1)\mod 3,0\leq k\leq2).\label{eq1.28}
\end{equation}
\end{theorem}

\begin{remark}\label{r1.2}
Analogously, for all $\lambda\in\mathbb{D}_{a/3}$,
$$v'_k(\lambda,x)-i\lambda\zeta_ke^{i\lambda\zeta_kx}\rightarrow0\quad(x\rightarrow\infty,0\leq k\leq2).$$
\end{remark}

For the solutions $\{u_k(\lambda,)\}_0^2$ to the boundary problem \eqref{eq1.13}, (b) \eqref{eq1.14}, analogously to \eqref{eq1.15}, the following integral equation is true:
\begin{equation}
u_k(\lambda,x)=e^{i\lambda\zeta_kx}+i\int\limits_{-\infty}^x\frac{s_t(i\lambda(x-t))}{(i\lambda)^2}q(t)u_k(\lambda,t)dt\quad(0\leq k\leq2)\label{eq1.29}
\end{equation}
which we can write as
$$(I-i\widehat{K})u_k(\lambda,x)=e^{i\lambda\zeta_kx}$$
where $\widehat{K}_\lambda$ is a Volterra operator,
$$(\widehat{K}_\lambda f)(x)\stackrel{\rm def}{=}\int\limits_{-\infty}^xK_1(\lambda,x,t)q(t)f(t)dt\quad(f\in L^2(\mathbb{R})$$
($K_1(\lambda,x,t)$ is given by \eqref{eq1.17}). Hence it follows that
$$u_k(\lambda,x)=\sum\limits_0^\infty i^n\widehat{K}_\lambda^ne^{i\lambda\zeta_kx},$$
and for the kernels $\widehat{K}_n(\lambda,x,t)$ estimates analogous to \eqref{eq1.21} are true, and
\begin{equation}
u_k(\lambda,x)=e^{i\lambda\zeta_nx}+\int\limits_{-\infty}^x\widehat{N}(\lambda,x,t)q(t)e^{i\lambda\zeta_kt}dt\quad(0\leq k\leq2)\label{eq1.30}
\end{equation}
where
$$\widehat{N}(\lambda,x,t)=\sum\limits_1^\infty i^n\widehat{K}_n(\lambda,x,t).$$

\begin{theorem}\label{t1.2}
Within the framework of assumptions of Theorem \ref{t1.1}, Jost solutions $\{u_k(\lambda,x)\}_0^2$ of the boundary value problem \eqref{eq1.13}, (b) \eqref{eq1.14} exist and are given by \eqref{eq1.30}. Functions $\{u_k(\lambda,x)\}_0^2$ are analytical with respect to $\lambda$ in the circle $\mathbb{D}_{a/3}$ \eqref{eq1.27} and
\begin{equation}
u_k(\lambda\zeta_1,x)=u_{k'}(\lambda,x)\quad(k'=(k+1)\mod3,0\leq k\leq2).\label{eq1.31}
\end{equation}
\end{theorem}

\begin{remark}\label{r1.3}
For all $\lambda\in\mathbb{D}_{a/3}$,
$$u'_k(\lambda,x)-i\lambda\zeta_ke^{i\lambda\zeta_kx}\rightarrow0\quad(x\rightarrow-\infty,0\leq k\leq2).$$
\end{remark}

\begin{lemma}\label{l1.3}
For all $\lambda\in\mathbb{D}_{a/3}\backslash\{0\}$, functions $\{v_k(\lambda,x)\}_0^2$ \eqref{eq1.24} ($\{u_k(\lambda,x)\}_0^2$ \eqref{eq1.30}) are linearly independent.
\end{lemma}

P r o o f. Lemma's statement follows from asymptotics \eqref{eq1.14}. Give a direct proof. Assuming contrary, suppose that for some $\lambda\in\mathbb{D}_{a/3}\backslash\{0\}$ there are numbers $\{\mu_k\}_0^2$ from $\mathbb{C}$ such that $\displaystyle{\sum\mu_kv_k(\lambda,x)=0}$ ($\forall x\in\mathbb{R}$), then
$$\sum\limits_k\mu_kv_k(\lambda,x)=0;\quad\sum\limits_k\mu_kv'_k(\lambda,x)=0;\quad\sum\limits_k\mu_kv''_k(\lambda,x)=0.$$
Determinant of this system
\begin{equation}
\Delta(\lambda,x)\stackrel{\rm def}{=}\det\left[
\begin{array}{cccc}
v_0(\lambda,x)&v_1(\lambda,x)&v_2(\lambda,x)\\
v'_0(\lambda,x)&v'_1(\lambda,x)&v'_2(\lambda,x)\\
v''_0(\lambda,x)&v''_1(\lambda,x)&v''_2(\lambda,x)
\end{array}\right]\label{eq1.32}
\end{equation}
does not depend on $x$ since $\Delta'(\lambda,x)=0$ ($\{v_k(\lambda,x)\}_0^2$ are solutions to \eqref{eq1.13}), therefore, in view of \eqref{eq1.14},
$$\Delta(\lambda,x)=\Delta_0(\lambda,x)+o\left(\frac1x\right)$$
where
$$\Delta_0(\lambda,x)=\det\left|
\begin{array}{ccc}
e^{i\lambda\zeta_0x}&e^{i\lambda\zeta_1x}&e^{i\lambda\zeta_2x}\\
i\lambda\zeta_0e^{i\lambda\zeta_0x}&i\lambda\zeta_1e^{i\lambda\zeta_1x}&i\lambda\zeta_2e^{i\lambda\zeta_2x}\\
(i\lambda\zeta_0)^2e^{i\lambda\zeta_0x}&(i\lambda\zeta_1)^2e^{i\lambda\zeta_1x}&(i\lambda\zeta_2)^2e^{i\lambda\zeta_2x}
\end{array}\right|=(i\lambda)^3\det\left|
\begin{array}{cccc}
1&1&1\\
1&\zeta_1&\zeta_2\\
1&\zeta_2&\zeta_1
\end{array}\right|$$
$$=-3\sqrt3\lambda^3.$$
So, $\Delta(\lambda,x)=-3\sqrt3\lambda^3$, therefore, $\Delta(\lambda,x)\not=0$ for all $\lambda\in\mathbb{D}_{a/3}\backslash\{0\}$, and thus $\mu_k=0$ ($\forall k$). $\blacksquare$
\vspace{5mm}

{\bf 1.3} Proceed to the analytical (by $\lambda$) properties of the solutions $\{v_k(\lambda,x)\}_0^2$ and $\{u_k(\lambda,x)\}_0^2$. Functions
\begin{equation}
\psi_k(\lambda,x)\stackrel{\rm def}{=}v_k(\lambda,x)e^{-i\lambda\zeta_kx}\quad(0\leq k\leq2)\label{eq1.33}
\end{equation}
due to \eqref{eq1.15}, are solutions to the integral equations
\begin{equation}
\psi_k(\lambda,x)=1-i\int\limits_x^\infty e^{i\lambda\zeta_k(t-x)}\frac{s_2(i\lambda(x-t))}{(i\lambda)^2}q(t)\psi_k(\lambda,t)dt\quad(0\leq k\leq2).\label{eq1.34}
\end{equation}
Kernel of this equation, for $k=0$, is
$$e^{i\lambda\zeta_0(t-x)}\frac{s_2(i\lambda(x-t))}{(i\lambda)^2}=\frac1{3(i\lambda)^2}\left\{1+\zeta_1e^{i\lambda(\zeta_0-\zeta_1)(t-x)}+\zeta_2e^{i\lambda(\zeta_0-\zeta_2)(t-x)}
\right\},$$
and since
\begin{equation}
\zeta_1-\zeta_0=-i\sqrt3\zeta_2;\quad\zeta_0-\zeta_2=-i\sqrt3\zeta_1;\quad\zeta_2-\zeta_1=-i\sqrt3\zeta_0,\label{eq1.35}
\end{equation}
and
$$i\lambda(\zeta_0-\zeta_1)=\frac{\sqrt3}2[(\alpha-\sqrt3\beta)+i(\alpha\sqrt3+\beta)];\quad i\lambda(\zeta_0-\zeta_2)=\frac{\sqrt3}2[(-\alpha-\sqrt3\beta)+i(\alpha\sqrt3+\beta)]$$
($\lambda=\alpha+i\beta$; $\alpha$, $\beta\in\mathbb{R}$), then
\begin{equation}
\left|e^{i\lambda(\zeta_0-\zeta_1)(t-x)}\right|\leq e^{\frac{\sqrt3}2(\alpha-\sqrt3\beta)(t-x)};\quad\left|e^{i\lambda(\zeta_0-\zeta_2)(t-x)}\right|\leq e^{-\frac{\sqrt3}2(\alpha+\beta\sqrt3)(t-x)}.\label{eq1.36}
\end{equation}

\begin{picture}(200,200)
\put(0,100){\vector(1,0){200}}
\put(100,0){\vector(0,1){200}}
\put(200,67){\line(-3,1){200}}
\put(0,67){\line(3,1){200}}
\put(120,150){$\Omega_0$}
\qbezier(153,120)(100,140)(49,120)
\qbezier(100,23)(155,35)(168,125)
\qbezier(45,120)(30,70)(100,25)
%\qbezier(125,100)(115,110)(120,104)
\put(140,105){$S_2(i)$}
\put(70,110){$S_4(i)$}
%\put(80,88){$S_4$}
\put(70,70){$S_6(i)$}
\put(50,100){$S_5(i)$}
\put(0,80){$\Omega_2$}
\put(115,70){$S_1(i)$}
\put(105,110){$S_3(i)$}
\put(180,80){$\Omega_1$}
%\put(165,50){$\Omega_1$}
%\put(80,40){$\psi_1^*(\lambda,0)$}
\end{picture}

\hspace{20mm} Fig. 1

So, for all $\lambda=\alpha+i\beta\in\mathbb{C}$ such that $\alpha-\beta\sqrt3<0$ and $-\alpha-\beta\sqrt3<0$, every exponent of the kernel is less than $1$ in modulo ($t-x>0$) and exponentially decreases as $t\rightarrow\infty$ when $\lambda$ belongs to this set. Applying the method of successive approximations to equation \eqref{eq1.34} ($k=0$), we obtain that $\psi_0(\lambda,x)$ is holomorphic in the domain $\alpha<\beta\sqrt3$, $-\alpha<\beta\sqrt3$. This set is described in terms of the sectors $\{S_p\}$ \eqref{eq1.9}. By $S_p(i)$, we denote sectors obtained from $S_p$ \eqref{eq1.9} via ``${\displaystyle-\frac\pi2}$'' rotation,
\begin{equation}
S_p(i)\stackrel{\rm def}{=}(-i)S_p\quad(1\leq p\leq6).\label{eq1.37}
\end{equation}
Define the sectors
\begin{equation}
\Omega_0\stackrel{\rm def}{=}S_3(i)\cup S_4(i)\cup(-i\widehat{l}_{\zeta_0});\quad\Omega_1\stackrel{\rm def}{=}S_1(i)\cup S_2(i)\cup(-i\widehat{l}_{\zeta_2});\quad\Omega_2\stackrel{\rm def}{=}S_5(i)\cup S_6(i)\cup(-i\widehat{l}_{\zeta_1}).\label{eq1.38}
\end{equation}
Domain of holomorphicity of $\psi_0(\lambda,x)$ ($\alpha<\beta\sqrt3$; $-\alpha<\beta\sqrt3$) coincides with $\Omega_0$.

Rewrite equation \eqref{eq1.34} ($k=0$) as
\begin{equation}
\psi_0(\lambda,x)=1+i\int\limits_x^\infty\frac1{3\lambda^2}\left[1+\zeta_1e^{i\lambda(\zeta_2-\zeta_1)(t-x)}+\zeta_2e^{i\lambda(\zeta_1-\zeta_2)(t-x)}\right]q(t)\psi_1(\lambda,t)
dt\label{eq1.39}
\end{equation}
and considering
\begin{equation}
\left|\frac1{3\lambda^2}\left(1+\zeta_1e^{i\lambda(\zeta_0-\zeta_1)(t-x)}+\zeta_2e^{i\lambda(\zeta_0-\zeta_2)(t-x)}\right)\right|<\frac1{\lambda^2}\quad(\forall\lambda\in\Omega_0)
\label{eq1.40}
\end{equation}
equation \eqref{eq1.39} implies that $|\psi_0(\lambda,x)|<(1-r(\lambda))^{-1}$, here $r(\lambda)=q_1/|\lambda|^2$, $q_1=\|q\|_{L^1}$ and $\lambda$ ($|\lambda|\gg1$) is chosen so that $r(\lambda)<1$. Therefore
\begin{equation}
|\psi_0(\lambda,x)-1|<r(\lambda)(1-r(\lambda))^{-1}\quad(\lambda\in\Omega_0,\,|\lambda|\gg1).\label{eq1.41}
\end{equation}
Relation \eqref{eq1.39} implies that
$$\psi_0(\lambda,x)=1+\frac i{3\lambda^2}\int\limits_x^\infty q(t)dt+C_0(\lambda)+D_0(\lambda)$$
where
\begin{equation}
\begin{array}{lll}
{\displaystyle C_0(\lambda,x)\stackrel{\rm def}{=}\frac i{3\lambda^2}\int\limits_x^\infty\left(1+\zeta_1e^{i\lambda(\zeta_0-\zeta_1)(t-x)}+\zeta_2e^{i\lambda(\zeta_0-\zeta_2)(t-x)}\right)q(t)(\psi_0(\lambda,t)-1)dt;}\\
{\displaystyle D_0(\lambda,x)\stackrel{\rm def}{=}\frac i{3\lambda^2}\int\limits_x^\infty\left(\zeta_1e^{i\lambda(\zeta_0-\zeta_1)(t-x)}+\zeta_2e^{i\lambda(\zeta_0-\zeta_2)(t-x)}\right)q(t)dt.}
\end{array}\label{eq1.42}
\end{equation}
Taking into account \eqref{eq1.42}, \eqref{eq1.41}, we obtain
$$|C_0(\lambda,x)|<r^2(\lambda)(1-r(\lambda))^{-1}\quad(\lambda\in\Omega_0,|\lambda\gg1).$$
For $D_0(\lambda,x)$ \eqref{eq1.42}, due to \eqref{eq1.36}, we have
\begin{equation}
|D_0(\lambda,x)|<\frac1{3|\lambda|^2}\left\{\int\limits_x^\infty e^{-\frac{\sqrt3}2(\beta\sqrt3-\alpha)(t-x)}|q(t)|dt-\int\limits_x^\infty e^{-\frac{\sqrt3}2(\beta\sqrt3+\alpha)(t-x)}|q(t)|dt\right\}.\label{eq1.43}
\end{equation}
The obvious inequality $1+kx\leq e^{kx}$ ($k$, $x\geq0$) implies $e^{-kx}\leq(1+kx)^{-1}$, and thus, for all $f\in L^2(\mathbb{R}_+)$, the following inequality holds:
$$\int\limits_0^\infty e^{-kx}|f(x)|dx\leq\int\limits_0^\infty\frac1{1+kx}|f(x)|dx\leq\left\{\int\limits_0^\infty|f(x)|^2dx\cdot\int\limits_0^\infty\frac{dx}{(1+kx)^2}\right\}^{1/2}=f_2\frac1{k^{1/2}}$$
($f_2=\|f\|_{L^2}$), and thus for \eqref{eq1.43}
$$|D_0(\lambda,x)|<\frac{2q_2}{3|\lambda|^2}\left\{\frac1{\sqrt3(\beta\sqrt3-|\alpha|)}\right\}^{1/2}=\frac{q_2}{|\lambda|^2}\left\{\frac8{9\sqrt3}\cdot\frac1{\sqrt3\beta-|\alpha|}
\right\}^{1/2}.$$

\begin{theorem}\label{t1.3}
For each $x\in\mathbb{R}$, solution $\psi_0(\lambda,x)$ to equation \eqref{eq1.34} ($k=0$) is analytical with respect to $\lambda$ in the sector $\Omega_0$ \eqref{eq1.38} and

(a) ${\displaystyle\psi_0(\lambda,x)=1+\frac i{3\lambda^2}\int\limits_x^\infty q(t)dt+C_0(\lambda,x)+D_0(\lambda,x)}$;
\begin{equation}
(b)\,v'_0(\lambda,x)e^{-i\lambda\zeta_0x}=i\lambda\zeta_0-\frac1{3\lambda}\int\limits_x^\infty q(t)dt+C_1(\lambda,x)+D_1(\lambda,x)\quad(\lambda\in\Omega_0);\label{eq1.44}
\end{equation}

(c) ${\displaystyle v''_0(\lambda,x)e^{-i\lambda\zeta_0x}=(i\lambda\zeta_0)^2-\frac i3\int\limits_x^\infty q(t)dt+C_2(\lambda,x)+D_2(\lambda,x).}$

For the functions $\{C_k(\lambda,x)\}_0^2$ and $\{D_k(\lambda,x)\}_0^2$ the following estimates hold:
\begin{equation}
|C_k(\lambda,x)|<|\lambda|^kr^2(\lambda)(1-r(\lambda))^{-1};\quad |D_k(\lambda,x)|<q_2|\lambda|^{k-2}\delta(\lambda)\quad(\lambda\in\Omega_2,|\lambda|\gg1)\label{eq1.45}
\end{equation}
($0\leq k\leq2$) where
\begin{equation}
r(\lambda)=\frac{q_1}{|\lambda|^2}\quad(q_s=\|q\|_{L^s}, s=1,2;|\lambda|\gg1;r(\lambda)<1);\quad\delta(\lambda)=(\beta\sqrt3-|\alpha|)^{-1/2}\label{eq1.46}
\end{equation}
($\lambda=\alpha+i\beta\in\Omega_0$).
\end{theorem}

\begin{remark}\label{r1.4}
Relation $\psi_k(\lambda,x)=\psi_0(\lambda_k,x)$ ($k=1$, $2$) implies that functions $\psi_1(\lambda,x)$ and $\psi_2(\lambda,x)$ are holomorphic in the sectors $\Omega_1$ and $\Omega_2$ \eqref{eq1.38} and analogues of the formulas \eqref{eq1.44} -- \eqref{eq1.46} are true in these sectors.

On the boundary of sector $\Omega_0$, $C_k(\lambda,x)\rightarrow0$ as $\lambda\rightarrow\infty$, and $D_k(\lambda,x)$ on $\partial\Omega_0$ also vanishes as $\lambda\rightarrow\infty$, due to Riemann -- Lebesgue lemma.

When $\lambda=i\omega$, ($\in\Omega_0$, $\omega>0$), function $C_0(\lambda,x)$ has order of smallness $|\omega|^{-4}$, and $D_0(\lambda,x)$, $|\omega|^{-2,5}$, correspondingly, therefore,
\begin{equation}
\lim\limits_{\omega\rightarrow\infty}3i\omega^2(\psi_0(i\omega,x)-1)=\int\limits_x^\infty q(t)dt.\label{eq1.47}
\end{equation}
\end{remark}
\vspace{5mm}

{\bf 1.4} Analogously to \eqref{eq1.33}, define the functions
\begin{equation}
\varphi_k(\lambda,x)\stackrel{\rm def}{=}u_k(\lambda,x)e^{-i\lambda_kx}\quad(0\leq k\leq2)\label{eq1.48}
\end{equation}
which, due to \eqref{eq1.29}, are solutions to the equations
\begin{equation}
\varphi_k(\lambda,x)=1+i\int\limits_{-\infty}^xe^{i\lambda\zeta_k(t-x)}\frac{s_2(i\lambda(x-t))}{(i\lambda)^2}q(t)\varphi_k(\lambda,t)dt\quad(0\leq k\leq2).\label{eq1.49}
\end{equation}
By $\Omega_k^-$, denote the centrally symmetric to $\Omega_k$ \eqref{eq1.18} sectors,
\begin{equation}
\Omega_k^-\stackrel{\rm def}{=}\{\lambda:-\lambda\in\Omega_k\}\quad(0\leq k\leq2).\label{eq1.50}
\end{equation}

\begin{theorem}\label{t1.4}
For all $x\in\mathbb{R}$, solution $\varphi_0(\lambda,x)$ to equation \eqref{eq1.49} ($k=0$) is holomorphic relative to $\lambda$ in sector $\Omega_0^-$ \eqref{eq1.50} and

(a) ${\displaystyle\varphi_0(\lambda,x)=1-\frac i{3\lambda^2}\int\limits_{-\infty}^xq(t)dt+\widehat{C}_0(\lambda,x)+\widehat{D}_(\lambda,x);}$
\begin{equation}
(b)\,u'_0(\lambda,x)e^{-i\lambda\zeta_0x}=i\lambda\zeta_0+\frac1{3\lambda}\int\limits_{-\infty}^xq(t)dt+\widehat{C}_1(\lambda,x)+\widehat{D}_1(\lambda,x)\quad(\lambda\in\Omega_0^-,
|\lambda|\gg1);\label{eq1.51}
\end{equation}

(c) ${\displaystyle u''_0(\lambda,x)e^{-i\lambda\zeta_0x}=(i\lambda\zeta_0)^2+\frac13\int\limits_{-\infty}^xq(t)dt+\widehat{C}_2(\lambda,x)+\widehat{D}_2(\lambda,x).}$\\
Functions $\{\widehat{C}_k(\lambda,x)\}_0^2$ and $\{\widehat{D}_k(\lambda,x)\}_0^2$ satisfy the estimates
\begin{equation}
|\widehat{C}_k(\lambda,x)|<|\lambda|^kr(\lambda)(1-r(\lambda))^{-1};\quad|\widehat{D}_k(\lambda,x)|<q_2|\lambda|^{k-2}\widehat{\delta}(\lambda)\quad(\lambda\in\Omega_0^{-1},|
\lambda|\gg1)\label{eq1.52}
\end{equation}
($0\leq k\leq2$) where $r(\lambda)$ is from \eqref{eq1.46} and $\widehat{\delta}(\lambda)=(-\beta\sqrt3-|\alpha|)^{-1/2}$ ($\lambda-\alpha+i\beta\in\Omega_0^-$).
\end{theorem}

\begin{remark}\label{r1.5}
Relation $\varphi_k(\lambda,x)=\varphi_0(\lambda\zeta_k,x)$ ($k=1$, $2$) implies that the functions $\varphi_1(\lambda,x)$ and $\varphi_2(\lambda,x)$ are analytical in $\Omega_1^-$ and $\Omega_2^-$, besides  formulas analogous to \eqref{eq1.51}, \eqref{eq1.52} are true.
\end{remark}

\section{Scattering problem}\label{s2}

{\bf 2.1} Write each Jost solution $\{u_k(\lambda,x)\}_0^3$ using the other Jost solutions $\{v_k(\lambda,x)\}_0^2$ (Lemma \ref{l1.2}),
\begin{equation}
u_k(\lambda,x)=\sum\limits_lt_{k,l}v_l(\lambda,x)\quad(0\leq k\leq2,\lambda\in\mathbb{D}_{a/3}\backslash\{0\})\label{eq2.1}
\end{equation}
or
\begin{equation}
u(\lambda,x)=T(\lambda)v(\lambda,x)\label{eq2.2}
\end{equation}
where $u(\lambda,x)\stackrel{\rm def}{=}\col[u_0(\lambda,x),u_1(\lambda,x),u_2(\lambda,x)]$; $v(\lambda,x)\stackrel{\rm def}{=}\col[v_0(\lambda,x),v_1(\lambda,x),$ $v_2(\lambda,x)]$, and $T(\lambda)$ is the {\bf transition matrix} \cite{1} -- \cite{3},
\begin{equation}
T(\lambda)\stackrel{\rm def}{=}\left[
\begin{array}{ccc}
t_{0,0}(\lambda)&t_{0,1}(\lambda)&t_{0,2}(\lambda)\\
t_{1,0}(\lambda)&t_{1,1}(\lambda)&t_{1,2}(\lambda)\\
t_{2,0}(\lambda)&t_{2,1}(\lambda)&t_{2,2}(\lambda)
\end{array}\right].\label{eq2.3}
\end{equation}

\begin{remark}\label{r2.1}
Equality
\begin{equation}
u_0(\lambda,x)=\sum\limits_lt_{0,l}(\lambda)v_l(\lambda,x)\label{eq2.4}
\end{equation}
(see \eqref{eq2.1}, $k=0$) implies the other $k=1$, $2$ relations \eqref{eq2.1}. Really, due to \eqref{eq1.28} and \eqref{eq1.31},
$$t_{1,k'}(\lambda)=t_{0,k}(\lambda\zeta_1)\,(k'=(k+1)\mod3);\quad t_{2,\widehat{k}}(\lambda)=t_{0,k}(\lambda\zeta_2)\,(\widehat{k}=(k+2)\mod3).$$
\end{remark}

Define the {\bf operation ``*''},
\begin{equation}
f^*(\lambda)\stackrel{\rm def}{=}\overline{f(\overline{\lambda})}.\label{eq2.5}
\end{equation}
Note that the operations ``*'' and $\lambda\rightarrow\lambda\zeta_1$ do not commute. By $W_{k,s}(v,\lambda,x)$, denote the Wronskian of functions $v_k(\lambda,x)$ and $v_s(\lambda,x)$:
\begin{equation}
W_{k,s}(v,\lambda,x)\stackrel{\rm def}{=}\{v_k(\lambda,x),v_s(\lambda,x)\}=v_k(\lambda,x)v'_s(\lambda,x)-v'_k(\lambda,x)v_s(\lambda,x)\,(0\leq k,s\leq2).\label{eq2.6}
\end{equation}

\begin{lemma}\label{l2.1}
Wronskians $\{W_{k,s}(v,\lambda,x)\}$ \eqref{eq2.6} have the representation
\begin{equation}
\begin{array}{ccc}
W_{0,1}(v,\lambda,x)=\sqrt3\lambda\zeta_2v_1^*(\lambda,x);\quad W_{1,2}(v,\lambda,x)=\sqrt3\lambda\zeta_0v_0^*(\lambda,x);\\ W_{2,0}(v,\lambda,x)=\sqrt3\lambda\zeta_1v_2^*(\lambda,x).
\end{array}\label{eq2.7}
\end{equation}
\end{lemma}

P r o o f. Since
$$W'_{k,s}(v,\lambda,x)=v_k(\lambda,x)v''_s(\lambda,x)-v''_k(\lambda,x)v_s(\lambda,x);$$
$$W''_{k,s}(v,\lambda,x)=v'_k(\lambda,x)v''_s(\lambda,x)-v''_k(\lambda,x)v'_s(\lambda,x);$$
then, upon differentiating the second equality ($\{v_k(\lambda,x)\}$ is the solution to \eqref{eq1.13}), we obtain that $W_{k,s}(v,\lambda,x)$ satisfies the equation
\begin{equation}
iy'''(x)-q(x)y(x)=-\lambda^3y(x)\label{eq2.8}
\end{equation}
which coincides with \eqref{eq1.13} upon substituting $\lambda\rightarrow\overline{\lambda}$ and taking complex conjugate. To obtain \eqref{eq2.7}, one has to take into account asymptotic (a) \eqref{eq1.14}. So, for $W_{0,1}(v,\lambda,x)$, as $x\rightarrow\infty$, we have $W_{0,1}(v,\lambda,x)\rightarrow i\lambda(\zeta_1-\zeta_2)e^{i\lambda(\zeta_0+\zeta_1)x}=\sqrt3\lambda\zeta_2e^{-i\lambda\zeta_2}$ (see \eqref{eq1.35}), which gives $W_{0,1}(v,\lambda,x)=\sqrt3\lambda\zeta_2v_1^*(\lambda,x)$. $\blacksquare$

\begin{remark}\label{r2.2}
For the Wronskians
\begin{equation}
W_{k,s}(u,\lambda,x)\stackrel{\rm def}{=}\{u_k(\lambda,x),u_s(\lambda,x)\}=u_k(\lambda,x)u'_s(\lambda,x)-u'_k(\lambda,x)u_s(\lambda,x)\,(0\leq k,s\leq2)\label{eq2.9}
\end{equation}
corresponding to the Jost solutions $\{u_k(\lambda,x)\}$, representations analogous to \eqref{eq2.7} are true.
\end{remark}

Using \eqref{eq2.1}, calculate the Wronskian
$$W_{1,2}(u,\lambda,x)=\{t_{1,0}(\lambda)v_0(\lambda,x)+t_{1,1}(\lambda)v_1(\lambda,x)+t_{1,2}(\lambda)v_2(\lambda,x),t_{2,0}(\lambda)v_0(\lambda,x)$$
$$+t_{2,1}(\lambda)v_1(\lambda,x)+t_{2,2}(\lambda)v_2(\lambda,x)\}=W_{0,1}(v,\lambda,x)\left|
\begin{array}{ccc}
t_{1,0}(\lambda)&t_{1,1}(\lambda)\\
t_{2,0}(\lambda)&t_{2,1}(\lambda)
\end{array}\right|$$
$$+W_{0,2}(v,\lambda,x)\left|
\begin{array}{cccc}
t_{1,0}(\lambda)&t_{1,2}(\lambda)\\
t_{2,0}(\lambda)&t_{2,2}(\lambda)
\end{array}\right|+W_{1,2}(v,\lambda,x)\left|
\begin{array}{ccc}
t_{1,1}(\lambda)&t_{1,2}(\lambda)\\
t_{2,1}(\lambda)&t_{2,2}(\lambda)
\end{array}\right|,$$
and, according to \eqref{eq2.7}, we have
\begin{equation}
\zeta_0u_0^*(\lambda,x)=\zeta_0v_0^*(\lambda,x)T_{0,0}(\lambda)+\zeta_2v_1^*(\lambda,x)T_{0,2}(\lambda)+\zeta_1v_2^*(\lambda,x)T_{0,1}(\lambda)\label{eq2.10}
\end{equation}
where $T_{k,s}(\lambda)$ is the algebraic complement of the elements $t_{k,s}(\lambda)$ of the matrix $T(\lambda)$ \eqref{eq2.3}. Analogously, considering $W_{0,1}(u,\lambda,x)$ and $W_{2,0}(u,\lambda,x)$, we obtain
\begin{equation}
\begin{array}{lll}
\zeta_2u_1^*(\lambda,x)=\zeta_0v_0^*(\lambda,x)T_{2,0}(\lambda)+\zeta_2v_1^*(\lambda,x)T_{2,2}(\lambda)+\zeta_1v_2^*(\lambda,x)T_{2,1}(\lambda);\\
\zeta_1u_2^*(\lambda,x)=\zeta_0v_0^*(\lambda,x)T_{1,0}(\lambda)+\zeta_2v_1^*(\lambda,x)T_{1,2}(\lambda)+\zeta_1v_2^*(\lambda,x)T_{1,1}(\lambda).
\end{array}\label{eq2.11}
\end{equation}
Note that relations \eqref{eq2.4}, \eqref{eq2.10} and \eqref{eq2.1}, \eqref{eq2.11} imply
\begin{equation}
\left\{
\begin{array}{lll}
t_{0,0}^*(\lambda)=T_{0,0}(\lambda);\quad t_{0,1}^*(\lambda)=\zeta_2T_{0,2};\quad t_{0,2}^*(\lambda)=\zeta_1T_{0,1};\\
t_{1,0}^*(\lambda)=\zeta_1T_{2,0}(\lambda);\quad t_{1,1}^*(\lambda)=T_{2,2}(\lambda);\quad t_{1,2}^*(\lambda)=\zeta_2T_{2,1}(\lambda);\\
t_{2,2}^*(\lambda)=\zeta_2T_{1,0}(\lambda);\quad t_{2,1}^*(\lambda)=\zeta_1T_{1,2}(\lambda);\quad t_{2,2}^*(\lambda)=T_{1,1}(\lambda).
\end{array}\right.\label{eq2.12}
\end{equation}
Equalities \eqref{eq2.10}, \eqref{eq2.11} in matrix form become
\begin{equation}
J^*U^*(\lambda,x)=\widehat{T}J^*v(\lambda,x)\label{eq2.13}
\end{equation}
where $J$ is the self-adjoint matrix
\begin{equation}
J=\left[
\begin{array}{cccc}
\zeta_0&0&0\\
0&0&\zeta_2\\
0&\zeta_1&0
\end{array}\right];\quad J^*=J;\quad J^2=I;\label{eq2.14}
\end{equation}
and $J^*$ is obtained from $J$ upon applying ``*'' to its elements, $\widehat{T}(\lambda)=[T_{k,l}(\lambda)]$ is the matrix formed by the algebraic complements to $T(\lambda)$ \eqref{eq2.3}. Using the fact that $\widehat{T}(\lambda)=\Delta(\lambda)(T^t(\lambda))^{-1}$ where $\Delta(\lambda)=\det T(\lambda)$ and $T^t(\lambda)$ is the matrix obtained from $T(\lambda)$ \eqref{eq2.3} via transposition, relation \eqref{eq2.13} yields
$$J^*T^t(\lambda)J^*u(\lambda,x)=\Delta(\lambda)v(\lambda,x),$$
or, upon applying ``*'',
$$\Delta^*(\lambda)v(\lambda,x)=J(T^t(\lambda))^*Ju(\lambda,x)=JT(\overline{\lambda})^+Ju(\lambda,x)$$
where $T^+(\lambda)$ is the matrix conjugate to $T(\lambda)$. Hence, upon finding $v(\lambda,x)$ and substituting it into \eqref{eq2.2}, we have $\{\Delta^*(\lambda)-T(\lambda)JT^+(\overline{\lambda})J\}u(\lambda,x)=0$, and thus
\begin{equation}
\Delta^*(\lambda)J=T(\lambda)JT^+(\overline{\lambda}).\label{eq2.15}
\end{equation}

\begin{lemma}\label{l2.2}
Matrix $T(\lambda)$ \eqref{eq2.3} is unimodular, $\det T(\lambda)=1$, and
\begin{equation}
J=T(\lambda)JT^+(\overline{\lambda})\quad(\forall\lambda\in\mathbb{D}_{a/3})\label{eq2.16}
\end{equation}
where $J$ is given by \eqref{eq2.14} and $\{\zeta_k\}_0^2$, correspondingly, by \eqref{eq1.2}. The following equality is true:
\begin{equation}
1=t_{0,0}(\lambda)t_{0,0}^*(\lambda)+\zeta_1t_{0,2}(\lambda)+J_2t_{0,1}(\lambda)t_{0,2}^*(\lambda)\quad(\forall\lambda\in\mathbb{D}_{a/3}).\label{eq2.17}
\end{equation}
\end{lemma}

P r o o f. Equation \eqref{eq2.15} implies
\begin{equation}
(\Delta^*(\lambda))^3=\Delta(\lambda)\cdot\Delta^*(\lambda).\label{eq2.18}
\end{equation}
Applying ``*'' to this equality and subtracting the obtained from this equality, we have $(\Delta^*(\lambda))^3=(\Delta(\lambda))^3$, and thus $\Delta^*=\zeta_k\Delta(\lambda)$ where $\zeta_k$ is one of roots of \eqref{eq1.2}. This fact and \eqref{eq2.18} imply that $\Delta(\lambda)=\zeta_k$. Now equality \eqref{eq2.15} becomes
$$\overline{\zeta}_kJ=T(\lambda)JT^*(\overline{\lambda}).$$
Upon equating the elements on the $(1,1)$ place in this equality, we have
$$\overline{\zeta}_k=t_{0,0}(\lambda)t_{0,0}^*(\lambda)+\zeta_1t_{0,2}(\lambda)t_{0,1}^*(\lambda)+\zeta_2t_{0,1}(\lambda)t_{0,2}^*(\lambda).$$
For $\lambda\in\mathbb{R}$, the right-hand part of this equality is real and thus $\zeta_k=\zeta_0$, i.e., matrix $T(\lambda)$ is unimodular, $\det T(\lambda)=1$. This fact and the last relation yield \eqref{eq2.17}. $\blacksquare$

Equating other matrix elements in equality \eqref{eq2.16}, we don't obtain any new relations, apart from \eqref{eq2.17}.

Relation \eqref{eq2.4} implies the equality
\begin{equation}
r_0(\lambda)u_0(\lambda,x)=v_0(\lambda,x)+sc_1(\lambda)v_1(\lambda,x)+sc_2(\lambda)v_2(\lambda,x)\label{eq2.19}
\end{equation}
where
\begin{equation}
r_0(\lambda)=t_{0,0}^{-1}(\lambda);\quad sc_1(\lambda)=t_{0,1}(\lambda)\cdot t_{0,0}^{-1}(\lambda);\quad sc_2(\lambda)=t_{0,2}(\lambda)t_{0,0}^{-1}(\lambda).\label{eq2.20}
\end{equation}

\begin{remark}\label{r2.3}
Asymptotic behavior of $r_0(\lambda)u_0(\lambda,x)$ \eqref{eq2.19} at ``$+\infty$'' is
$$r_0(\lambda)u_0(\lambda,x)\rightarrow e^{i\lambda\zeta_0x}+sc_1e^{i\lambda\zeta_1x}+sc_2e^{i\lambda\zeta_2x}\quad(x\rightarrow+\infty,\lambda\in\mathbb{D}_{a/3})$$
and thus the {\bf wave $e^{i\lambda\zeta_0x}$ falling from the right}  has the {\bf reflected wave} which is represented by the {\bf span} of waves $sc_1(\lambda)e^{i\lambda\zeta_1x}+sc_2(\lambda)e^{i\lambda\zeta_2x}$ where $sc_1(\lambda)$ and $sc_2(\lambda)$ are {\bf scattering coefficients of the wave} $e^{i\lambda\zeta_0x}$. At ``$-\infty$'', function $r_0(\lambda)u_0(\lambda,x)$ is
$$r_0(\lambda)u_0(\lambda,x)\rightarrow r_0(\lambda)e^{i\lambda\zeta_0x}\quad(x\rightarrow-\infty,\lambda\in\mathbb{D}_{a/3})$$
and describes the {\bf passed wave} where $r(\lambda)$ is the {\bf transition coefficient of the wave} $e^{i\lambda\zeta_2x}$.

{\bf Unitarity condition} \cite{1} -- \cite{3} of such scattering problem follows from \eqref{eq2.17}:
\begin{equation}
r_0(\lambda)r_0^*(\lambda)=1+\zeta_1sc_2(\lambda)sc_1^*(\lambda)+\zeta_2sc_1(\lambda)sc_2^*(\lambda)\quad(\lambda\in\mathbb{D}_{a/3}).\label{eq2.21}
\end{equation}

Analogous relations for falling waves $e^{i\lambda\zeta_1x}$ and $e^{i\lambda\zeta_2x}$ follow from \eqref{eq2.19} upon substituting $\lambda\rightarrow\lambda\zeta_1$ and $\lambda\rightarrow\lambda\zeta_2$.
\end{remark}
\vspace{5mm}

{\bf 2.2} Calculate the coefficients $\{t_{0,k}(\lambda)\}_0^2$ of decomposition \eqref{eq2.4}. Consider the system
$$\left\{
\begin{array}{lll}
u_0(\lambda,x)=t_{0,0}(\lambda)v_0(\lambda,x)+t_{0,1}(\lambda)v_1(\lambda,x)+t_{0,2}(\lambda)v_2(\lambda,x);\\
u'_0(\lambda,x)=t_{0,0}(\lambda)v'_0(\lambda,x)+t_{0,1}(\lambda)v'_1(\lambda,x)+t_{0,2}(\lambda)v'_2(\lambda,x);\\
u''_0(\lambda,x)=t_{0,0}(\lambda)v''_0(\lambda,x)+t_{0,1}(\lambda)v''_1(\lambda,x)+t_{0,2}(\lambda)v''_2(\lambda,x).
\end{array}\right.$$
Its determinant $\Delta(\lambda,x)$ \eqref{eq1.38} equals $\Delta(\lambda,x)=-3\sqrt3\lambda^3$, and determinant $\Delta_0(\lambda,x)$ corresponding to $t_{0,0}(\lambda)$ is
\begin{equation}
\begin{array}{ccc}
{\displaystyle\Delta_0(\lambda,x)\stackrel{\rm def}{=}\det\left[
\begin{array}{ccc}
u_0(\lambda,x)&v_1(\lambda,x)&v_2(\lambda,x)\\
u'_0(\lambda,x)&v'_1(\lambda,x)&v'_2(\lambda,x)\\
u''_0(\lambda,x)&v''_1(\lambda,x)&v''_2(\lambda,x)
\end{array}\right]=u_0(\lambda,x)W''_{1,2}(v,\lambda,x)}\\
-u'_0(\lambda,x)W'_{1,2}(v,\lambda,x)+u''_0(\lambda,x)W_{1,2}(v,\lambda,x);
\end{array}\label{eq2.22}
\end{equation}
and, taking \eqref{eq2.7} into account, we have
\begin{equation}
\Delta_0(\lambda,x)=\sqrt3\lambda\zeta_0\{u_0(\lambda,x)(v_0^*(\lambda,x))''-u'_0(\lambda,x)(v_0^*(\lambda,x))'+u''_0(\lambda,x)v_0^*(\lambda,x)\}.\label{eq2.23}
\end{equation}
Function $\Delta_0(\lambda,x)$ \eqref{eq2.23} does not depend on $x$ since $\Delta'_0(\lambda,x)=0$ because $u_0(\lambda,x)$ and $v_0^*(\lambda,x)$ are solutions to the equations \eqref{eq1.13} and \eqref{eq2.8}. Hereinafter, we consider $x\equiv0$.

\begin{lemma}\label{l2.3}
Coefficients $\{t_{0,k}(\lambda)\}_0^2$ of decomposition \eqref{eq2.4} are
\begin{equation}
t_{0,k}(\lambda)=-\frac{\zeta_k}{3\lambda^2}\{u_0(\lambda,0)(v_{n(k)}^*(\lambda,0))''-u'_0(\lambda,0)(v_{n(k)}^*(\lambda,0))'+u''_0(\lambda,0)v_{n(k)}^*(\lambda,0)\}\label{eq2.24}
\end{equation}
where $n(0)=0$, $n(1)=2$, $n(2)=1$ ($0\leq k\leq2$).
\end{lemma}

Proceed to the set of zeros of the function $t_{0,0}(\lambda)$. Theorems \ref{t1.3}, \ref{t1.4} imply that the function
\begin{equation}
t_{0,0}(\lambda)=-\frac1{3\lambda^2}\{u_0(\lambda,0)(v_0^*(\lambda,0))''-u'_0(\lambda,0)(v_0^*(\lambda,0))'+u''_0(\lambda,0)v_0^*(\lambda,0)\}\label{eq2.25}
\end{equation}
is holomorphic in the sector $\Omega_0^-$ \eqref{eq1.50}. Let $z\in\Omega_0^-$ be a zero of the function $t_{0,0}(\lambda)$ \eqref{eq2.25}, $t_{0,0}(z)=0$, then $z$ also is a zero of the function $\Delta_0(\lambda,x)$ \eqref{eq2.23}, $\Delta_0(z,x)=0$. Linear independence of the second and the third columns of matrix \eqref{eq2.23} and $\Delta_0(z,x)=0$ imply that the first column is a linear span of these columns. and thus
\begin{equation}
u_0(z,x)=c_1(z)v_1(z,x)+c_2(z)v_2(z,x).\label{eq2.26}
\end{equation}
Note that \eqref{eq2.26} coincides with \eqref{eq2.4} since $t_{0,0}(z)=0$.

\begin{picture}(200,200)
\put(0,100){\vector(1,0){200}}
\put(100,0){\vector(0,1){200}}
\put(100,100){\vector(1,-2){50}}
\put(100,100){\vector(-1,-2){50}}
\put(100,100){\vector(3,-1){100}}
\put(100,100){\vector(-3,-1){100}}
\put(80,190){$il_{\zeta_0}$}
\put(190,110){$l_{\zeta_0}$}
\put(2,50){$il_{\zeta_1}$}
\put(71,50){$\times$}
\put(75,44){$\varkappa\zeta_2$}
\put(105,30){$\Omega_0^-$}
\put(120,50){$\circ$}
\put(124,55){$-\widehat{\varkappa}\zeta_2$}
\put(56,3){$l_{\zeta_2}$}
\put(157,5){$\widehat{l}_{\zeta_1}$}
\put(190,70){$il_{\zeta_2}$}
\qbezier(153,80)(100,60)(49,80)
\end{picture}

\hspace{20mm} Fig. 2

Function $u_0(z,x)$ is an eigenfunction of the operator $L_q$ \eqref{eq1.10}, $L_qu_0(z,x)=z^3u_0(z,x)$ if only $u_0(z,x)\in L^2(\mathbb{R})$. Hence we conclude that $z^3$ is real and thus either $z=\varkappa\zeta_2\in l_{\zeta_2}$ ($\varkappa>0$, $z^3>0$) or $z=-\widehat{\varkappa}\zeta_1\in\widehat{l}_{\zeta_1}$ ($\widehat{\varkappa}>0$, $z^3<0$). It is left to make sure that $u_0(z,x)\in L^2(\mathbb{R})$. When $x\rightarrow-\infty$, function $u_0(z,x)$, while exponentially decreasing, tends to zero (see Section \ref{s1}), therefore $\chi_-u_0(z,x)\in L^2(\mathbb{R}_-)$ ($\chi_\pm$ are characteristic functions of $\mathbb{R}_\pm$). The right-hand side of equality \eqref{eq2.26} is holomorphic in the disk $\mathbb{D}_{a/3}$. Function $c_1(z)v_1(z,x)$ in \eqref{eq2.26}, for $z=\varkappa\zeta_2\in l_{\zeta_2}$ has an oscillating summand $c_1(z)e^{i\varkappa x}$, and function $c_2(z)v_2(z,x)$ for these $z$ ($=\varkappa\zeta_2$) vanishes exponentially decreasing. Therefore, in order to $u_0(z,x)\in L^2(\mathbb{R}_+)$ for $z\in l_{\zeta_2}$, it is necessary that $c_1(z)=t_{1,1}(z)=0$ and thus $u_0(z,x)=t_{0,2}(z)v_2(z,x)$ ($c_2(z)=t_{2,2}(z)$). Analogously, inclusion $u_0(z,x)\in L^2(\mathbb{R}_+)$ for $z\in\widehat{l}_{\zeta_1}$ holds if $c_2(z)=0$ and $u_0(z,x)=t_{0,1}(z)v_1(z,x)$.

\begin{lemma}\label{l2.4}
If \eqref{eq1.12} is true and
\begin{equation}
3\cdot\max\{2\sqrt{q_1},2\sqrt{q_2},9\}<a\quad(q_1\stackrel{\rm def}{=}\|q\|_{L^1},\,q_2\stackrel{\rm def}{=}\|q\|_{L^2}),\label{eq2.27}
\end{equation}
then the set of zeros of the function $t_{0,0}(\lambda)$ \eqref{eq2.25} is finite and is given by $\Lambda\stackrel{\rm def}{=}\Lambda_2\cup\widehat{\Lambda}_1$ ($\subset\Omega_0^-\cap\mathbb{D}_{a/3}$ and
\begin{equation}
\begin{array}{cccc}
\Lambda_2\stackrel{\rm def}{=}\{\varkappa_l\zeta_2\in l_{\zeta_2}:0<\varkappa_l<a/3;1\leq m\in\mathbb{Z}_+\};\\
\widehat\Lambda_1\stackrel{\rm def}{=}\{-\widehat{\varkappa_s}\zeta_1\in\widehat{l}_{\zeta_1}:0<\widehat{\varkappa}_s<a/3;1\leq s\leq\widehat{m}\in\mathbb{Z}_+\}.
\end{array}\label{eq2.28}
\end{equation}
Zeros of function $t_{0,0}(\lambda)$ are simple. Function $u_0(z,x)\in L^2(\mathbb{R})$ is an eigenfunction for the operator $L_q$ \eqref{eq1.10}, besides, if $z\in\Lambda_2$, then $z$ are joint zeros of $t_{0,0}(\lambda)$ and $t_{0,2}(\lambda)$, and $u_0(z,x)=t_{0,1}(z)v_1(z,x)$.
\end{lemma}

P r o o f. First, we show that each summand in representation \eqref{eq2.25} of the function $t_{0,0}(\lambda)$ is given by ${\displaystyle\frac13+\delta(\lambda)}$ where ${\displaystyle|\delta(\lambda)|<\frac13}$, when $\lambda\in l_{\zeta_2}\cup\widehat{l}_{\zeta_1}$ and $|\lambda|\gg1$. Consider the first summand in \eqref{eq2.25}
$$A_1(\lambda)\stackrel{\rm def}{=}-\frac1{3\lambda^2}u_0(\lambda,0)(v_0^*(\lambda,0))''.$$
Equation \eqref{eq1.51} implies that
\begin{equation}
u_0(\lambda,0)=1-\frac i{3\lambda^2}\int\limits_{-\infty}^0q(t)dt+M(\lambda),\label{eq2.29}
\end{equation}
and, according to \eqref{eq1.52},
$$M(\lambda)\stackrel{\rm def}{=}\widehat{C}_0(\lambda,0)+\widehat{D}_0(\lambda,0);\quad|M(\lambda)|<r^2(\lambda)(1-r(\lambda))^{-1}+\widehat{\delta}(\lambda)q_2|\lambda|^{-2}\quad(|\lambda|\gg1).$$
Let $q_0\stackrel{\rm def}{=}\max\{q_1,q_2\}$, then for all $\lambda$ ($|\lambda|\gg1$), $r_0(\lambda)<q_0/|\lambda|^2$ ($<1$) and
$$|M(\lambda)|<r_0(\lambda)(r_0(\lambda)/(1-r_0(\lambda))+\widehat{\delta}(\lambda)).$$
Analogously, due to \eqref{eq1.44}, \eqref{eq1.45}, we obtain
\begin{equation}
(v_0^*(\lambda,0))''=-\lambda^2\left(1-\frac i{3\lambda^2}\int\limits_0^\infty q(t)dt+M_2^*(\lambda)\right)\label{eq2.30}
\end{equation}
where $M_2^*(\lambda)=-(C_2^*(\lambda,0)+D_2^*(\lambda,0))\lambda^{-2}$ and
$$|M_2^*(\lambda)|<r_0(\lambda)(r_0(\lambda)/(1-r_0(\lambda))+\delta^*(\lambda)).$$
Substituting \eqref{eq2.29}, \eqref{eq2.30} into representation of $A_1(\lambda)$, we have
\begin{equation}
A_1(\lambda)=\frac13(1+\widehat{M}(\lambda))\label{eq2.31}
\end{equation}
where
$$\widehat{M}(\lambda)=-\frac i{3\lambda^2}\int\limits_{-\infty}^\infty q(t)dt-\left(\frac i{3x^2}\right)^2\int\limits_0^\infty q(t)dt\int\limits_{-\infty}^0q(t)dt+M(\lambda)+M_2^*(\lambda)$$
$$-\frac i{3\lambda^2}\int\limits_{-\infty}^0q(t)dtM_2^*(\lambda)-\frac i{3\lambda^2}\int\limits_0^\infty q(t)dt\cdot M(\lambda)+M(\lambda)M_2^*(\lambda).$$
If $|M(\lambda)|<m(\lambda)$ and $|M_2^*(\lambda)|<m(\lambda)$, then we obtain that
\begin{equation}
 |\widehat{M}(\lambda)|<\frac13r_0(\lambda)+\frac19r_0^2(\lambda)+2m(\lambda)+\frac13r_0(\lambda)m(\lambda)+m^2(\lambda).\label{eq2.32}
\end{equation}
Let $\lambda=-\tau-i\sqrt3\tau\in l_{\zeta_2}$ ($\tau>0$), then $|\lambda|=2\tau$ and $\delta(\lambda)=\delta^*(\lambda)=(\sqrt{2\tau})^{-1}$, therefore
$$r_0(\lambda)=\frac{q_0}{(2\tau)^2};\quad\delta(\lambda)=\delta^*(\lambda)=\frac1{\sqrt{2\tau}};\quad m(\lambda)=r_0(\lambda)\left(\frac{r_0(\lambda)}{1-r_0(\lambda)}+\frac1{\sqrt{2\tau}}\right).$$
Choose $\lambda$ (or $\tau>0$) so that ${\displaystyle r_0(\lambda)<\frac14}$ and ${\displaystyle\frac1{\sqrt{2\tau}}<\frac12}$, then ${\displaystyle m(\lambda)<\frac16}$, and \eqref{eq2.32} yields that
$$|\widehat{M}(\lambda)|<\frac14\left(\frac13+\frac1{36}+\frac43+\frac1{18}+\frac19\right)<1.$$
Consequently, for such $\lambda$ (or $\tau$) on the ray $l_{\zeta_2}$, function $A_1(\lambda)$ is ${\displaystyle A_1(\lambda)=\frac13(1+\delta(\lambda))}$ and $|\delta(\lambda)|<1$. Conditions ${\displaystyle r_0(\lambda)<\frac14}$ and ${\displaystyle\frac1{\sqrt{2\tau}}<\frac13}$ give $2\sqrt{q_0}<2\tau$ and $9<2\tau$, and thus $\max\{2\sqrt{q_0},9\}<2\tau$, and since ${\displaystyle|\lambda|=2\tau<\frac a3}$, then hence follows condition \eqref{eq2.27}. For the second and the third summands in \eqref{eq2.25}, considerations are analogous. On the ray $\widehat{l}_{\zeta_1}$, finiteness of zero of the function $t_{0,0}(\lambda)$, under condition \eqref{eq2.27}, is proven in exactly the same way.

Prove that zeros of $t_{0,0}(\lambda)$ are simple. Let $z\in\Lambda$ be a zero of the function $t_{0,0}(\lambda)$, $t_{0,0}(z)=0$. Equality $L_qu_0(z,x)=z^3u_0(z,x)$, upon differentiation by $z$, implies
$$(L_q-z^3I)\dot u_0(z,x)=3z^2u_0(z,x)\quad\left(\dot u_0(z,x)=\frac d{dz}u_0(z,x)\right).$$
Scalar multiplying this equality by $u_0(z,x)$, we have
$$3z^2\|u_0(z,x)\|_{L^2}^2=\langle(L_0-z^3I)\dot u_0(z,x),u_0(z,x)\rangle,$$
and, upon integrating by parts,
$$3z^2\|u_0(z,x)\|_{L^2}^2=\left.i\{\dot u''_0(z,x)\overline{u_0(z,x)}-\dot u'_0\overline{u'_0(z,x)}+\dot u_0(z,x)\overline{u''_0(z,x)}\}\right|_{-\infty}^\infty.$$
At the lower limit ``$-\infty$'', expression at the right-hand side of this equality vanishes since $u_0(z,x)$ tends to zero while exponentially decreasing as $x\rightarrow-\infty$. Let $z=\varkappa\zeta_2\in\Lambda_2$, then \eqref{eq2.4}, taking into account $t_{0,0}(z)=t_{0,1}(z)=0$, gives
$$u_0(z,x)=t_{0,2}(z)v_2(z,x);\quad\dot u_0(z,x)=\dot t_{0,0}(z)v_0(z,x)+\dot t_{0,1}(z)v_1(z,x)$$
$$+\dot t_{0,2}(z)v_2(z,x)+t_{0,2}(z)\dot v_2(z,x),$$
use the asymptotic formulas \eqref{eq1.44}
$$v_0(z,x)=e^{i\varkappa\zeta_2x}+o(1);\quad v_1(z,x)=e^{i\varkappa x}+o(1);\quad v_2(z)=e^{i\varkappa\zeta_1x}+o(1)\quad(x\rightarrow\infty).$$
Then, after elementary transformations, we have that
$$\dot u''_0(z,x)\overline{u_0(z,x)}-\dot u'_0(z,x)\overline{u'_0(z,x)}+\dot u_0(z,x)\overline{u''_0(z,x)}=-3\varkappa^2\zeta_1\dot t_{0,0}(z)\overline{t_{0,2}}(z)+o(1)$$
(as $x\rightarrow\infty$), therefore
$$3\varkappa^2\zeta_1\|u_0(z,x)\|_{L^2}^2=-i^2\varkappa^2\zeta_1\dot t_{0,1}(z)\overline{t_{0,2}}(z).$$
Hence it follows that $\dot t_{0,0}(z)\not=0$. For $z=-\widehat{\varkappa}\zeta_1\in\widehat{\Lambda}_1$, proof is analogous.

\begin{remark}\label{r2.4}
Search for exact constants in condition \eqref{eq2.27} is not essential.
\end{remark}

\begin{remark}\label{r2.5}
If \eqref{eq1.12} and \eqref{eq2.27} are true, then operator $L_q$ \eqref{eq1.10} has a finite number of {\bf bound states} which are enumerated via points of the sets $\Lambda_2$ and $\widehat{\Lambda}_1$ \eqref{eq2.28}.

Function $t_{0,1}(\lambda)$ ($t_{0,2}(\lambda)$) \eqref{eq2.24} is analytical in the disk $\mathbb{D}_{a/3}$ and is defined only on the ray $il_{\zeta_1}$ ($il_{\zeta_2}$) dividing the sectors $\Omega_0^-$ and $\Omega_1^-$ ($\Omega_0^-$ and $\Omega_2^-$).
\end{remark}
\vspace{5mm}

{\bf 2.3} Equality \eqref{eq2.19} generates jump problems for functions holomorphic in adjoint sectors of a half-plane. So, \eqref{eq2.19} implies
$$r_0(\lambda)\{u_0(\lambda,x),v_2(\lambda,x)\}=W_{0,2}(v,\lambda,x)+sc_1W_{1,2}(v,\lambda,x)$$
and, taking into account \eqref{eq2.7},
$$r_0(\lambda)\{u_0(\lambda,x),v_2(\lambda,x)\}=-\sqrt3\lambda\zeta_1v_2^*(\lambda,x)+sc_1(\lambda)\sqrt3\lambda\zeta_0v_0^*(\lambda,x),$$
or, in terms of $\{\psi_k(\lambda,x)\}$ \eqref{eq1.33},
\begin{equation}
\psi_2^*(\lambda,\lambda)-r_0(\lambda)f_{0,2}(\lambda,x)=s_1(\lambda,x)\psi_0^*(\lambda,x)\quad(\lambda\in il_{\zeta_1});\label{eq2.33}
\end{equation}
where
\begin{equation}
s_1(\lambda,x)\stackrel{\rm def}{=}\zeta_2e^{i\lambda(\zeta_1-\zeta_0)x}sc_1(\lambda);\quad f_{0,2}(\lambda,x)\stackrel{\rm def}{=}-\frac{\zeta_2}{\sqrt3\lambda}e^{i\lambda\zeta_1x}\{u_0(\lambda,x),v_2(\lambda,x)\}.\label{eq2.34}
\end{equation}

\begin{picture}(400,200)
\put(24,100){\vector(1,0){200}}
\put(124,0){\vector(0,1){200}}
\put(124,100){\vector(1,-2){50}}
\put(124,100){\vector(-1,-2){50}}
\put(124,100){\vector(3,-1){100}}
\put(124,100){\vector(-3,-1){100}}
\put(26,50){$il_{\zeta_1}$}
\put(89,40){$\times$}
\put(99,44){$\varkappa\zeta_2$}
\put(144,50){$\circ$}
\put(150,50){$-\widehat{\varkappa}\zeta_2$}
\put(160,35){$\Omega_0^-\cap\Omega_1$}
\put(90,3){$l_{\zeta_2}$}
\put(181,15){$r_0(\lambda)f_{0,1}(\lambda,x)$}
\put(214,70){$il_{\zeta_2}$}
\qbezier(124,60)(159,50)(184,80)
\qbezier(124,63)(89,50)(64,80)
\qbezier(124,170)(74,140)(61,78)
\qbezier(124,170)(174,140)(187,78)
\put(88,130){$\Omega_1^-$}
\put(134,130){$\Omega_2^-$}
\put(174,130){$\psi_1^*(\lambda,x)$}
\put(24,130){$\psi_2^*(\lambda,x)$}
\put(49,50){$\Omega_0^-\cap\Omega_2$}
\put(0,20){$r_0(\lambda)f_{0,2}(\lambda,x)$}
\end{picture}

\hspace{20mm} Fig. 3

Equality \eqref{eq2.33} is a jump problem on the ray $il_{\zeta_1}$ for holomorphic at adjoint sectors $\Omega_1^-$ and $\Omega_0^-\cap\Omega_2$ of the left half-plane functions $\psi_2^*(\lambda,x)$ and $r_0(\lambda)f_{0,2}(\lambda,x)$. Function $\psi_2^*(\lambda,x)$ is holomorphic at the sector $\Omega_1^-$ and tends to $1$ when $\lambda\rightarrow\infty$ ($\lambda\in\Omega_1^-$), due to \eqref{eq1.14} and Remark \ref{r1.4}. Theorems \ref{t1.3}, \ref{t1.4} imply that the function $\{u_0(\lambda,x),v_2(\lambda,x)\}$ is analytical in the sector $\Omega_0^-\cap\Omega_2$. Show that $f_{0,2}(\lambda,x)$ \eqref{eq2.34} also tends to $1$ when $\lambda\rightarrow\infty$ ($\lambda\in\Omega_0^-\cap\Omega_2$). Since
$$f_{0,2}(\lambda,x)=-\frac{\zeta_2}{\sqrt3\lambda}e^{-i\lambda(\zeta_0+\zeta_2)x}\{u_0(\lambda,x),v_2(\lambda,x)\}$$
$$=-\frac{\zeta_2}{\sqrt3\lambda}\{\varphi_0(\lambda,x)v'_2(\lambda,x)
e^{-i\lambda\zeta_2x}-u'_0(\lambda,x)e^{-i\lambda\zeta_0x}\cdot\psi_2(\lambda,x)\},$$
then, taking into account \eqref{eq1.44} and \eqref{eq1.51}, we obtain
$$f_{0,2}(\lambda,x)=-\frac{\zeta_2}{\sqrt3\lambda}\{(1+o(\lambda^{-2}))(i\lambda\zeta_2+o(\lambda^{-1})-(i\lambda\zeta_0+o(\lambda^{-1}))(1+o(\lambda^{-2})))\}$$
$$=-\frac{\zeta_2}{\sqrt3\lambda}\{i\lambda(\zeta_2-\zeta_0)+o(\lambda^{-1})\}=1+o(\lambda^{-2}),$$
and thus $f_{0,2}(\lambda,x)\rightarrow1$ when $\lambda\rightarrow\infty$. Using $t_{0,0}(\lambda)\rightarrow1$ ($\lambda\rightarrow\infty$, $\lambda\in\Omega_0^-$), we have that $r_0(\lambda)\rightarrow1$ when $\lambda\rightarrow\infty$ ($\lambda\in\Omega_0^-$).

\begin{lemma}\label{l2.5}
For holomorphic at the sectors $\Omega_1^-$ and $\Omega_0^-\cap\Omega_2$ functions $\psi_0^*(\lambda,x)$ and $r_0(\lambda)f_{0,2}(\lambda,x)$ (excluding poles at the points of $\Lambda_2$ \eqref{eq2.28}), jump problem \eqref{eq2.33} takes place. These functions tend to $1$ when $\lambda\rightarrow\infty$ (inside the sectors $\Omega_1^-$ and $\Omega_0^-\cap\Omega_2$).
\end{lemma}

Analogously,
$$r_0(\lambda)\{u_0(\lambda,x),v_1(\lambda,x)\}=W_{0,1}(v,\lambda,x)+sc_2W_{2,1}(v,\lambda,x),$$
and thus
\begin{equation}
\psi_1^*(\lambda,x)-r_0(\lambda)f_{0,1}(\lambda,x)=s_2(\lambda,x)\psi_0^*(\lambda,x)\quad(\lambda\in il_{\zeta_2})\label{eq2.35}
\end{equation}
where
\begin{equation}
s_2(\lambda,x)\stackrel{\rm def}{=}\zeta_1e^{i\lambda(\zeta_2-\zeta_0)x}sc_2(\lambda);\quad f_{0,1}(\lambda,x)\stackrel{\rm def}{=}\frac{\zeta_1}{\sqrt3\lambda}e^{i\lambda\zeta_2x}\{u_0(\lambda,x),v_1(\lambda,x)\}.\label{eq2.36}
\end{equation}
As a result, we have the jump problem \eqref{eq2.35} on the ray $il_{\zeta_2}$ for holomorphic in the sectors $\Omega_2^-$ and $\Omega_0^-\cap\Omega_1$ of the right half-plane functions $\psi_1^*(\lambda,x)$ and $r_0(\lambda)f_{0,1}(\lambda,x)$ (excluding poles at the points of $\widehat{\Lambda}_1$ \eqref{eq2.28}. These functions tend to $1$ as $\lambda\rightarrow\infty$ (inside of $\Omega_2^-$ and $\Omega_0^-\cap\Omega_1$).

Again, from \eqref{eq2.19} we obtain
$$r_0(\lambda)\{u_0(\lambda,x),v_0(\lambda,x)\}=sc_1(\lambda)W_{1,0}(v,\lambda,x)+sc_2(\lambda)W_{2,0}(v,\lambda,x);$$
i.e.,
\begin{equation}
s_2(\lambda,x)\psi_2^*(\lambda,x)-s_1(\lambda,x)\psi_1^*(\lambda,x)=r_0(\lambda)f_{0,0}(\lambda,x),\label{eq2.37}
\end{equation}
here $s_1(\lambda,x)$ and $s_2(\lambda,x)$ are from \eqref{eq2.34} and \eqref{eq2.36}, and
\begin{equation}
f_{0,0}(\lambda,x)\stackrel{\rm def}{=}\frac{\zeta_0}{\sqrt3\lambda}e^{-2i\lambda\zeta_0x}\{u_0(\lambda,x),v_0(\lambda,x)\}\label{eq2.38}
\end{equation}
Unlike $f_{0,2}(\lambda,x)$ and $f_{0,1}(\lambda,x)$, function $f_{0,0}(\lambda,x)$ does not have analytic continuation into the sector.
\vspace{5mm}

{\bf 2.4.} Substituting $\lambda\rightarrow\lambda\zeta_1$ in \eqref{eq2.33}, we obtain
\begin{equation}
\psi_1^*(\lambda,x)-r_0(\lambda\zeta_1)f_{1,0}(\lambda,x)=s_1(\lambda\zeta_1,x)\psi_2^*(\lambda,x)\quad(\lambda\in il_{\zeta_0})\label{eq2.39}
\end{equation}
where
\begin{equation}
s_1(\lambda\zeta_1,x)\stackrel{\rm def}{=}\zeta_2e^{i\lambda(\zeta_2-\zeta_1)x}sc_1(\lambda\zeta_1);\quad f_{1,0}(\lambda,x)\stackrel{\rm def}{=}-\frac{\zeta_1}{\sqrt3\lambda}e^{i\lambda\zeta_2x}\{u_1(\lambda,x),v_0(\lambda,x)\}.\label{eq2.40}
\end{equation}

\begin{picture}(300,200)
\put(70,100){\vector(1,0){200}}
\put(170,0){\vector(0,1){200}}
\put(170,100){\vector(-1,2){50}}
\put(70,133){\vector(3,-1){200}}
\put(170,100){\vector(-3,-1){100}}
\put(100,190){$l_{\zeta_1}$}
\put(185,150){$\psi_1^*(\lambda,x)$}
\put(172,190){$il_{\zeta_0}$}
\put(140,147){$\times$}
\put(147,147){$\varkappa\zeta_1$}
\put(190,50){$\psi_0^*(\lambda,x)$}
\put(76,85){$\Omega_2\cap\Omega_1^-$}
\put(72,57){$il_{\zeta_1}$}
\put(0,105){$r_0(\lambda\zeta_1)f_{1,2}(\lambda,x)$}
\put(120,96){$\circ$}
\put(104,105){$-\widehat\varkappa\zeta_0$}
\put(104,125){$\Omega_0\cap\Omega_1^-$}
\put(247,58){$il_{\zeta_2}$}
\put(30,150){$r_0(\lambda\zeta_1)f_{1,0}(\lambda,x)$}
\qbezier(193,89)(206,140)(170,133)
\qbezier(203,90)(170,50)(139,88)
\qbezier(149,108)(130,105)(139,88)
\qbezier(150,107)(150,120)(170,130)
\put(177,110){$\Omega_2^-$}
\put(174,80){$\Omega_0^-$}
\put(50,130){$i\widehat{l}_{\zeta_2}$}
\end{picture}

\hspace{20mm} Fig. 4.

Equality \eqref{eq2.39} gives a jump problem on the ray $il_{\zeta_0}$ for holomorphic in the sectors $\Omega_2^-$ and $\Omega_0\cap\Omega_1^-$ functions $\psi_1^*(\lambda,x)$ and $r_0(\lambda\zeta_1)f_{1,0}(\lambda,x)$ (excluding poles in $\Lambda_1\stackrel{\rm def}{=}\zeta_2\Lambda_2$) which tend to $1$ when $\lambda\rightarrow\infty$ inside these sectors.

Analogously, upon substituting $\lambda\rightarrow\zeta_1\lambda$ in \eqref{eq2.35}, we obtain
\begin{equation}
\psi_0^*(\lambda,x)-r_0(\lambda\zeta_1)f_{1,2}(\lambda,x)=s_2(\lambda\zeta_1,x)\psi_2^*(\lambda,x)\quad(\lambda\in il_{\zeta_1})\label{eq2.41}
\end{equation}
where
\begin{equation}
s_2(\lambda\zeta_1,x)\stackrel{\rm def}{=}\zeta_1e^{i\lambda(\zeta_0-\zeta_1)x}sc_2(\lambda\zeta_1);\quad f_{1,2}(\lambda,x)\stackrel{\rm def}{=}\frac{\zeta_0}{\sqrt3\lambda}e^{i\lambda\zeta_0x}\{u_1(\lambda,x),v_2(\lambda,x)\}.\label{eq2.42}
\end{equation}
Relation \eqref{eq2.41} is a jump problem on the ray $il_{\zeta_1}$ for holomorphic in sectors $\Omega_0^-$ and $\Omega_2\cap\Omega_1^-$ functions $\psi_0^*(\lambda,x)$ and $r(\lambda\zeta_1)f_{1,2}(\lambda,x)$ (excluding poles in $\widehat{\Lambda}_0\stackrel{\rm def}{=}\zeta_2\widehat{\Lambda}_1$) and they tend to $1$ inside these sectors.

Substitution $\lambda\rightarrow\lambda\zeta_2$ in \eqref{eq2.33} gives us
\begin{equation}
\psi_0^*(\lambda,x)-r_0(\lambda\zeta_2)f_{2,1}(\lambda,x)=s_1(\lambda\zeta_2,x)\psi_1^*(\lambda,x)\quad(\lambda\in il_{\zeta_2})\label{eq2.43}
\end{equation}
where
\begin{equation}
s_1(\lambda\zeta_2,x)\stackrel{\rm def}{=}\zeta_2e^{i\lambda(\zeta_0-\zeta_2)x}sc_1(\lambda\zeta_2);\quad f_{2,1}(\lambda,x)\stackrel{\rm def}{=}-\frac{\zeta_0}{\sqrt3\lambda}e^{i\lambda\zeta_0x}\{u_2(\lambda,x),v_2(\lambda,x)\}.\label{eq2.44}
\end{equation}

\begin{picture}(300,200)
\put(70,100){\vector(1,0){200}}
\put(170,0){\vector(0,1){200}}
\put(170,100){\vector(1,2){50}}
\put(260,131){\vector(-3,-1){200}}
\put(170,100){\vector(3,-1){100}}
\put(168,170){$-\widehat\varkappa\zeta_2$}
\put(172,190){$il_{\zeta_0}$}
\put(198,162){$\circ$}
\put(210,155){$r_0(\lambda\zeta_2)f_{2,0}(\lambda,x)$}
\put(200,140){$\Omega_0\cap\Omega_2^-$}
\put(190,50){$\psi_0^*(\lambda,x)$}
\put(148,83){$\Omega_1^-$}
\put(148,103){$\Omega_1^-$}
\put(72,57){$il_{\zeta_1}$}
\put(224,104){$\varkappa\zeta_0$}
\put(260,120){$r_0(\lambda\zeta_2)f_{2,1}(\lambda,x)$}
\put(260,80){$\Omega_1\cap\Omega_2^-$}
\put(220,96){$\times$}
\put(104,125){$\psi_2^*(\lambda,x)$}
\put(247,58){$il_{\zeta_2}$}
\qbezier(193,110)(206,140)(170,133)
\qbezier(195,110)(220,120)(202,91)
\qbezier(203,88)(170,50)(139,88)
\qbezier(135,93)(130,120)(170,130)
\end{picture}

\hspace{20mm} Fig. 5.

We have the jump problem on the ray $il_{\zeta_2}$ for holomorphic in the sectors $\Omega_0^-$ and $\Omega_1\cap\Omega_2^-$ functions $\psi_0^*(\lambda,x)$ and $r_0(\lambda\zeta_2)f_{2,1}(\lambda,x)$ (excluding poles at the points of $\Lambda_0\stackrel{\rm def}{=}\zeta_1\Lambda_2$), and these functions tend to $1$ inside these sectors. Similar, upon substituting $\lambda\rightarrow\lambda\zeta_2$ in \eqref{eq2.35}, we have
\begin{equation}
\psi_2^*(\lambda,x)-r_0(\lambda\zeta_2)f_{2,0}(\lambda,x)=s_2(\lambda\zeta_2,x)\psi_1^*(\lambda,x)\quad(\lambda\in il_{\zeta_0})\label{eq2.45}
\end{equation}
where
\begin{equation}
s_2(\lambda\zeta_2)\stackrel{\rm def}{=}\zeta_1e^{i\lambda(\zeta_1-\zeta_2)x}sc_2(\lambda\zeta_2);\quad f_{2,0}(\lambda,x)=\frac{\zeta_2}{\sqrt3\lambda}e^{i\lambda\zeta_1x}\{u_2(\lambda,x),v_0(\lambda,x)\}.\label{eq2.46}
\end{equation}
So, we have the jump problem on the ray $il_{\zeta_0}$ for the analytic in the sectors $\Omega_1^-$ and $\Omega_0\cap\Omega_2^-$ functions $\psi_2^*(\lambda,x)$ and $r_0(\lambda\zeta_2)f_{2,0}(\lambda,x)$ (excluding poles at the points of $\widehat{\Lambda}_2=\zeta_1\widehat{\Lambda}_1$) which tend to $1$ as $\lambda\rightarrow\infty$ (inside these sectors).
\vspace{5mm}

{\bf 2.5} Consider the function
\begin{equation}
\mathcal{F}(\lambda,x)\stackrel{\rm def}{=}\left\{
\begin{array}{lllll}
\psi_0^*(\lambda,x)\quad(\lambda\in\Omega_0^-);\\
r_0(\lambda\zeta_2)f_{2,1}(\lambda,x)\quad(\lambda\in\Omega_1\cap\Omega_2^-);\\
r_0(\lambda\zeta_2)f_{2,0}(\lambda,x)\quad(\lambda\in\Omega_0\cap\Omega_2^-);\\
r_0(\lambda\zeta_1)f_{1,0}(\lambda,x)\quad(\lambda\in\Omega_0\cap\Omega_1^-);\\
r_0(\lambda\zeta_1)f_{1,2}(\lambda,x)\quad(\lambda\in\Omega_2\cap\Omega_1^-).
\end{array}\right.\label{eq2.47}
\end{equation}
on the rays $il_{\zeta_1}$ and $il_{\zeta_2}$, due to \eqref{eq2.41} and \eqref{eq2.43}, for the components of function $\mathcal{F}(\lambda,x)$ \eqref{eq2.47}, we have the jump problems
\begin{equation}
\left[
\begin{array}{lll}
\psi_0^*(\lambda,x)-r_0(\lambda\zeta_1)f_{1,2}(\lambda,x)=s_2(\lambda\zeta_1,x)\psi_2^*(\lambda,x)\quad(\lambda\in il_{\zeta_1});\\
\psi_0^*(\lambda,x)-r_0(\lambda\zeta_2)f_{2,1}(\lambda,x)=s_1(\lambda\zeta_2,x)\psi_1^*(\lambda,x)\quad(\lambda\in il_{\zeta_2}).
\end{array}\right.\label{eq2.48}
\end{equation}
Functions $r_0(\lambda\zeta_1)f_{1,0}(\lambda,x)$ and $r_0(\lambda\zeta_1)f_{1,2}(\lambda,x)$ are meromorphic in the sectors $\Omega_0\cap\Omega_1^-$ and $\Omega_2\cap\Omega_1^-$ and on its boundary (ray $i\widehat{l}_{\zeta_2}$) have the following boundary conditions:
\begin{equation}
\left[
\begin{array}{lll}
\left.r_0(\lambda\zeta_1)f_{1,0}(\lambda,x)\right|_{\lambda\in i\widehat{l}_{\zeta_2}}=r_0(-i\tau)f_{0,2}(-i\tau,x);\\
\left.r_0(\lambda\zeta_1)f_{1,2}(\lambda,x)\right|_{\lambda\in i\widehat{l}_{\zeta_2}}=r_0(-i\tau)f_{0,1}(-i\tau,x);
\end{array}\right.\label{eq2.49}
\end{equation}
where $\lambda=-i\tau\zeta_2$ ($\tau>0$). Analogously, for the meromorphic in sectors $\Omega_1\cap\Omega_2^-$ and $\Omega_0\cap\Omega_2^-$ functions $r_0(\lambda\zeta_2)f_{2,1}(\lambda,x)$ and $r_0(\lambda\zeta_2)f_{2,0}(\lambda,x)$, the boundary values on the ray $i\widehat{l}_{\zeta_1}$ are
\begin{equation}
\left[
\begin{array}{lll}
\left.r_0(\lambda\zeta_2)f_{2,1}(\lambda,x)\right|_{\lambda\in i\widehat{l}_{\zeta_1}}=r_0(-i\tau)f_{0,2}(-i\tau,x);\\
\left.r_0(\lambda\zeta_2)f_{2,0}(\lambda,x)\right|_{\lambda\in i\widehat{l}_{\zeta_1}}=r_0(-i\tau)f_{0,1}(-i\tau,x)
\end{array}\right.\label{eq2.50}
\end{equation}
($\lambda=-i\tau\zeta_1\in i\widehat{l}_{\zeta_1}$). Hence it follows that
\begin{equation}
\left[
\begin{array}{lll}
\left.r_0(\lambda\zeta_1)f_{1,2}(\lambda,x)\right|_{\lambda\in i\widehat{l}_{\zeta_2}}=\left.r_0(\lambda\zeta_2)f_{2,0}(\lambda,x)\right|_{\lambda\in i\widehat{l}_{\zeta_1}};\\
\left.r_0(\lambda\zeta_2)f_{2,1}(\lambda,x)\right|_{\lambda\in i\widehat{l}_{\zeta_2}}=\left.r_0(\lambda\zeta_1)f_{1,0}(\lambda,x)\right|_{\lambda\in i\widehat{l}_{\zeta_2}}.
\end{array}\right.\label{eq2.51}
\end{equation}

Each point $\lambda\in\Omega_0\cap\Omega_2^-$ we match with the symmetric relative to straight line $il_{\zeta_0}$ point $\lambda^+$ ($=-\overline\lambda$) from sector $\Omega_0\cap\Omega_1^-$. This involution $\lambda\rightarrow\lambda^+$ leaves ray $il_{\zeta_0}$ fixed. A holomorphic function $f(\lambda)$ in $\Omega_0\cap\Omega_2^-$ we match with the analytic function $f^+(\lambda)$ in $\Omega_0\cap\Omega_1^-$,
\begin{equation}
f^+(\lambda)=f(\lambda^+)\quad(\lambda^+\in\Omega_0\cap\Omega_2^-).\label{eq2.52}
\end{equation}
Relation \eqref{eq2.51} yields that the function
\begin{equation}
\mathcal{F}^+(\lambda,x)\stackrel{\rm def}{=}\left\{
\begin{array}{lllll}
\psi_0^*(\lambda,x)&(\lambda\in\Omega_0^-);\\
r_0(\lambda\zeta_2)f_{2,1}(\lambda,x)&(\lambda\in\Omega_1\cap\Omega_2^-);\\
(r_0(\lambda\zeta_1)f_{1,0}(\lambda,x))^+&(\lambda\in\Omega_0\cap\Omega_2^-);\\
(r_0(\lambda\zeta_2)f_{2,0}(\lambda,x))^+&(\lambda\in\Omega_0\cap\Omega_1^-);\\
r_0(\lambda\zeta_1)f_{1,2}(\lambda,x)&(\lambda\in\Omega_2\cap\Omega_1^-)
\end{array}\right.\label{eq2.53}
\end{equation}
is holomorphic in $\Omega_0^-$ and meromorphic in $\Omega_1^-$ ($\Omega_2^-$) excluding poles at $\{-\widehat{\varkappa}_s\zeta_0\}_1^{\widehat m}$, $\{-\widehat\varkappa_s\zeta_1\}_1^{\widehat m}$ ($\{\varkappa_l\zeta_0\}_1^m$, $\{\varkappa_l\zeta_2\}_1^m$). On the rays $il_{\zeta_1}$ and $il_{\zeta_2}$, the jump problem \eqref{eq2.48} holds, and on the ray $il_{\zeta_0}$, according to \eqref{eq2.39}, \eqref{eq2.45}, such problem is given by
\begin{equation}
r_0(\lambda\zeta_2)f_{2,0}(\lambda,x)-r_0(\lambda\zeta_1)f_{1,0}(\lambda,x)=(s_1(\lambda\zeta_1,x)-1)\psi_2^*(\lambda,x)-(s_2(\lambda\zeta_2,x)-1)\psi_1^*(\lambda,x)\label{eq2.54}
\end{equation}
($\lambda\in il_{\zeta_0}$).

\begin{picture}(300,200)
\put(70,100){\vector(1,0){200}}
\put(170,0){\vector(0,1){200}}
\put(170,100){\vector(1,2){50}}
\put(170,100){\vector(-1,2){50}}
\put(260,131){\vector(-3,-1){200}}
\put(70,133){\vector(3,-1){200}}
\put(168,170){$-\varkappa_l\zeta_2$}
\put(172,190){$il_{\zeta_0}$}
\put(198,162){$\times$}
\put(0,110){$r_0(\lambda\zeta_1)f_{1,2}(\lambda,x)$}
\put(200,140){$(r_0(\lambda\zeta_1)f_{1,0}(\lambda,x))^+$}
\put(190,50){$\psi_0^*(\lambda,x)$}
\put(148,83){$\Omega_0^-$}
\put(80,97){$\circ$}
\put(72,57){$il_{\zeta_1}$}
\put(224,104){$\varkappa_l\zeta_0$}
\put(260,120){$r_0(\lambda\zeta_2)f_{2,1}(\lambda,x)$}
\put(220,96){$\times$}
\put(133,162){$\circ$}
\put(140,155){$\widehat\varkappa_s\zeta_1$}
\put(10,150){$(r_0(\lambda\zeta_2)f_{2,0}(\lambda,x))^+$}
\put(247,58){$il_{\zeta_2}$}
\qbezier(193,110)(206,140)(170,133)
\qbezier(195,110)(220,120)(196,91)
\qbezier(203,88)(170,50)(139,88)
\qbezier(115,83)(95,100)(120,114)
\qbezier(120,114)(140,145)(170,130)
\put(70,87){$-\widehat{\varkappa}_s\zeta_0$}
\end{picture}

\hspace{20mm} Fig. 6 .

Since $\mathcal{F}^+(\lambda,x)\rightarrow1$ when $\lambda\rightarrow\infty$ (inside $\{\Omega_k^-\}_0^2$),
$$\mathcal{F}^+(\lambda,x)=1+\sum\limits_1^m\frac{R_l(x)}{\lambda-\varkappa_l}+\sum\limits_1^m\frac{R'_l(x)}{\lambda+\varkappa_l\zeta_2}+\sum\limits_1^{\widehat{m}}
\frac{\widehat{R}_s(x)}{\lambda+\widehat\varkappa_s}+\sum\limits_1^{\widehat m}\frac{\widehat{R'}_s(x)}{\lambda-\widehat\varkappa_s\zeta_1}$$
\begin{equation}
+\frac1{2\pi i}\int\limits_{il_{\zeta_1}}\frac{s_1(\lambda\zeta_1,x)\psi_2^*(z,x)}{z-\lambda}dz-\frac1{2\pi i}\int\limits_{il_{\zeta_2}}\frac{s_1(\lambda\zeta_2,x)\psi_1^*(\lambda,x)}{z-\lambda}dz\label{eq2.55}
\end{equation}
$$+\frac1{2\pi i}\int\limits_{il_{\zeta_0}}\frac{(s_1(\lambda\zeta_1,x)-1)\psi_2^*(z,x)-(s_2(z\zeta_1,x)-1)\psi_1^*(z,x)}{z-\lambda}dz$$
where $R_l(x)$, $R'_l(x)$ ($\widehat{R_s}(x)$, $\widehat{R'}_s(x)$) are residues of the function $\mathcal{F}^+(\lambda,x)$ \eqref{eq2.53} at the points $\lambda=\varkappa_l$, $\lambda=-\varkappa_l\zeta_2$ ($\lambda=-\widehat\varkappa_s$, $\lambda=\widehat\varkappa_s\zeta_1$). Signs `$-$' and `$+$' on the sides of rays $\{il_{\zeta_k}\}_0^2$ go counterclockwise as well as the signs on $l_{\zeta_0}$ ($=R_+$) in accordance with half-planes $\mathbb{C}_-$ and $\mathbb{C}_+$. Taking into account that the function $r_0(\lambda\zeta_2)f_{2,0}(\lambda,x)$ is obtained from $r_0(\lambda\zeta_1)f_{1,2}(\lambda,x)$ via the substitution $\lambda\rightarrow\lambda\zeta_2$ and operation `$+$' \eqref{eq2.52} changes the sign before residue, then $\widehat{R}'_s(x)=-R_s(x)$ ($R'_l(x)=-R_l(x)$), therefore
$$\mathcal{F}^+(\lambda,x)=1-\zeta_1\sum\limits_1^m\frac{\varkappa_lR_l(x)}{(\lambda-\varkappa_l)(\lambda+\varkappa_l\zeta_2)}+\zeta_2\sum\limits_1^{\widehat m}\frac{\widehat\varkappa_s\widehat{R}_s(x)}{(\lambda+\widehat\varkappa_s)(\lambda-\widehat\varkappa_s\zeta_1)}$$
\begin{equation}
+\frac1{2\pi i}\int\limits_0^\infty s_2(i\zeta_2\tau,x)\psi_1^*(i\tau,x)\frac{d\tau}{\tau+i\zeta_2\lambda}-\frac1{2\pi i}\int\limits_0^\infty s_1(i\zeta_1\tau,x)\psi_2^*(i\tau,x)\frac{d\tau}{\tau+i\zeta_1\lambda}\label{eq2.56}
\end{equation}
$$+\frac1{2\pi i}\int\limits_0^\infty\{(s_1(i\zeta_1\tau,x)-1)\psi_2^*(i\tau,x)-(s_2(i\zeta_2\tau,x)-1)\psi_1^*(i\tau,x)\}\frac{d\tau}{\tau+i\lambda}.$$
Use the identity
$$\int\limits_0^\infty f(x)\frac{dx}{x+z}=\int\limits_0^\infty f(x)\int\limits_x^\infty\frac{dt}{(t+z)^2}dx=\int\limits_0^\infty\int\limits_0^tf(x)dx\frac{dt}{(t+z)^2}\quad(z\not\in\mathbb{R}_+),$$
then \eqref{eq2.56}, for $\lambda\in\Omega_0^-$, implies
$$\psi_0^*(\lambda,x)=1-\zeta_1\sum\limits_1^m\frac{\varkappa_lR_l(x)}{(\lambda-\varkappa_l)(\lambda+\varkappa_l\zeta_2)}+\zeta_2\sum\limits_1^{\widehat m}\frac{\widehat\varkappa_s\widehat{R}_s(x)}{(\lambda+\widehat\varkappa_s)(\lambda-\widehat\varkappa_s\zeta_1)}$$
$$+\frac1{2\pi i}\int\limits_0^\infty\int\limits_0^\tau s_2(i\zeta_2t,x)\psi_1^*(it,x)dt\frac{d\tau}{(\tau+i\zeta_2\lambda)^2}$$
\begin{equation}
-\frac1{2\pi i}\int\limits_0^\infty\int\limits_0^\tau s_1(i\zeta_1t,x)\psi_2^*(it,x)dt\frac{d\tau}{(\tau+i\zeta_1\lambda)^2}\label{eq2.57}
\end{equation}
$$+\frac1{2\pi i}\int\limits_0^\infty\int\limits_0^\tau s_1(i\zeta_1t,x)\psi_2^*(it,x)dt\frac{d\tau}{(\tau+i\lambda)^2}-\frac1{2\pi i}\int\limits_0^\infty\int\limits_0^\tau s_2(i\zeta_2t,x)\psi_1^*(it,x)dt\frac{d\tau}{(\tau+i\lambda)^2}$$
$$+\frac1{2\pi i}\int\limits_0^\infty\int\limits_0^\tau(\psi_1^*(it,x)-\psi_2^*(it,x))\frac{d\tau}{(\tau+i\lambda)^2}.$$
Calculate boundary values of the function $\psi_1^*(\lambda,x)$ ($=\mathcal{F}^+(\lambda,x)$ \eqref{eq2.56} as $\lambda\in\Omega_0^-$) when $\lambda\rightarrow i\zeta_1t$ ($\in il_{\zeta_1}$) and when $\lambda\rightarrow i\zeta_2t$ ($\in il_{\zeta_2}$) assuming that $\lambda\in\Omega_0^-$,
\begin{equation}
\left\{
\begin{array}{lll}
{\displaystyle\psi_2^*(it,x)=1+\zeta_2\sum\limits_1^m\frac{\varkappa_lR_l(x)}{(t+i\zeta_2\varkappa_l)(t-i\zeta_1\varkappa_l)}-\zeta_0\sum\limits_1^{\widehat m}\frac{\widehat\varkappa_s\widehat{R}_s(x)}{(t-i\zeta_2\widehat\varkappa_s)(t+i\widehat\varkappa_s)}}\\
{\displaystyle+\frac1{2\pi i}\int\limits_0^\infty s_2(i\zeta_2\tau,x)\psi_1^*(i\tau,x)\frac{d\tau}{\tau-(t+i0)}-\frac1{2\pi i}\int\limits_0^\infty s_1(i\zeta_1\tau,x)\psi_2^*(i\tau,x)\frac{d\tau}{\tau-\zeta_2t}}\\
{\displaystyle+\frac1{2\pi i}\int\limits_0^\infty\{(s_1(i\zeta_1\tau,x)-1)\psi_2^*(i\tau,x)-(s_2(i\zeta_2\tau,x)-1)\psi_1^*(i\tau,x)\}\frac{d\tau}{\tau-\zeta-1t};}\\
{\displaystyle\psi_1^*(it,x)=1+\zeta_0\sum\limits_1^m\frac{\varkappa_lR_l(x)}{(t+i\zeta_1\varkappa_l)(t-i\varkappa_l)}-\zeta_1\sum\limits_1^{\widehat m}\frac{\widehat\varkappa_s\widehat{R}_s(x)}{(t-i\zeta_1\widehat\varkappa_s)(t+i\zeta_2\widehat\varkappa_s)}}\\
{\displaystyle+\frac1{2\pi i}\int\limits_0^\infty s_2(i\zeta_2\tau,x)\frac{d\tau}{\tau-\zeta_1t}-\frac1{2\pi i}\int\limits_0^\infty s_1(i\zeta_1\tau,x)\psi_2^*(i\tau,x)\frac{d\tau}{\tau-(t-i0)}}\\
{\displaystyle+\frac1{2\pi i}\int\limits_0^\infty\{(s_1(i\zeta_1\tau,x)-1)\psi_2^*(i\tau,x)-(s_2(i\zeta_2\tau,x)-1)\psi_1^*(i\tau,x)\}\frac{d\tau}{\tau-\zeta_2t};}
\end{array}\right.\label{eq2.58}
\end{equation}
where $t\pm i0$ corresponds to the boundary values on $\mathbb{R}_+$ from the half-planes $\mathbb{C}_\pm$.

Residue of the function $r_0(\zeta_2)f_{2,1}(\lambda,x)$ at the point $\lambda=\varkappa_l$ coincides with the residue of the function $r_0(\lambda)f_{0,2}(\lambda,x)$ at the point $\lambda=\varkappa_l\zeta_2$,
$$R_l(x)=\res\limits_{\varkappa_l}\{r_0(\lambda\zeta_2)f_{2,1}(\lambda,x)\}=\res\limits_{\varkappa_l\zeta_2}\{r_0(\lambda)f_{0,2}(\lambda,x)\},$$
and, taking into account \eqref{eq2.33} and meromorphicity of $r_0(\lambda)f_{0,2}(\lambda,x)$ in the sector $\Omega_0^-\cap\Omega_2$, we obtain
$$R_l(x)=-\res\limits_{\varkappa_l\zeta_2}\{s_1(\lambda,x)\psi_0^*(\lambda,x)\}=-\res\limits_{\varkappa_l\zeta_2}\left\{\zeta_2e^{i\lambda(\zeta_1-\zeta_0)x}\frac{t_{0,1}(\lambda)
\cdot(\lambda-\varkappa_l\zeta_2)}{t_{0,0}(\lambda)}\frac{\mathcal{F}^+(\lambda,x)}{\lambda-\varkappa_l\zeta_2}\right\}$$
(see \eqref{eq2.20}, \eqref{eq2.34}). So,
$$R_l(x)=-\zeta_2e^{i\varkappa_l(\zeta_0-\zeta_2)x}\frac{t_{0,1}(\varkappa_l\zeta-2)}{t'_{0,0}(\varkappa_l\zeta_2)}\mathcal{F}^+(\varkappa_l\zeta_2,x),$$
and thus
\begin{equation}
\mathcal{F}^+(\varkappa_l\zeta_2,x)=b_le^{i\varkappa_l(\zeta_2-\zeta_0)x}R_l(x)\quad\left(b_l\stackrel{\rm def}{=}-\zeta_1\frac{t'_{0,0}(\varkappa_l\zeta_2)}{t_{0,1}(\varkappa_l\zeta_2)}\right).\label{eq2.59}
\end{equation}
Similarly, it is proved that
\begin{equation}
\mathcal{F}^+(-\widehat\varkappa_s\zeta_1,x)=\widehat{b}_se^{-i\widehat\varkappa_s(\zeta_1-\zeta_0)}\widehat{R}_s(x)\quad\left(\widehat{b}_s\stackrel{\rm def}{=}-\zeta_2\frac{t'_{0,0}(-\widehat\varkappa_s\zeta_1)}{t_{0,2}(-\widehat{\varkappa}_s\zeta_1)}\right).\label{eq2.60}
\end{equation}
Equalities \eqref{eq2.59}, \eqref{eq2.60} imply another $m+\widehat{m}$ equations:
\begin{equation}
\left\{
\begin{array}{llllll}
{\displaystyle b_le^{i\varkappa_l(\zeta_2-\zeta_0)x}R_l(x)=1+\zeta_2\sum\limits_{p=1}^m\frac{\varkappa_pR_p(x)}{(\varkappa_p-\varkappa_l\zeta_2)(\varkappa_p+\varkappa_l)}}\\
{\displaystyle-\zeta_1
\sum\limits_1^{\widehat m}\frac{\widehat\varkappa_s\widehat{R}_s(x)}{(\widehat\varkappa_s+\varkappa_l\zeta_2)(\widehat\varkappa_s-\varkappa_l\zeta_1)}}
{\displaystyle+\frac1{2\pi i}\int\limits_0^\infty s_2(i\zeta_2\tau,x)\psi_1^*(i\tau,x)\frac{d\tau}{\tau+i\zeta_1\varkappa_l}}\\
{\displaystyle-\frac1{2\pi i}\int\limits_0^\infty s_1(i\zeta_1\tau,x)\psi_2^*(i\tau,x)\frac{d\tau}{\tau+i\varkappa_l}}\\
{\displaystyle+\frac1{2\pi i}\int\limits_0^\infty\{(s_1(i\zeta_1\tau,x)-1)\psi_2^*(i\tau,x)-(s_2(i\zeta_2\tau,x)-1)\psi_1^*(i\tau,x)\}\frac{d\tau}{\tau+i\zeta_2\varkappa_l}}\\
(1\leq l\leq m);\\
{\displaystyle\widehat{b}_se^{-i\widehat\varkappa_s(\zeta_1-\zeta_0)x}\widehat{R}_s(x)=1+\zeta_2\sum\limits_1^m\frac{\varkappa_lR_l(x)}{(\varkappa_l+\widehat\varkappa_s\zeta_1)
(\varkappa_l-\widehat\varkappa_s\zeta_2)}}\\
{\displaystyle+\zeta_1\sum\limits_{p=1}^{\widehat m}\frac{\widehat{\varkappa}_p\widehat{R}_p(x)}{(\widehat{\varkappa}_p-\widehat{\varkappa}_s\zeta_1)(\widehat{\varkappa}_p+\widehat{\varkappa}_s)}}\\
{\displaystyle+\frac1{2\pi i}\int\limits_0^ls_2(i\zeta_2\tau,x)\psi_1^*(i\tau,x)\frac{d\tau}{\tau-i\widehat{\varkappa}_s}-\frac1{2\pi i}\int\limits_0^\infty s_1(i\zeta_1\tau,x)\psi_2^*(i\tau,x)\frac{d\tau}{\tau-i\zeta_2\widehat{\varkappa}_s}}\\
{\displaystyle+\frac1{2\pi i}\int\limits_0^\infty\{(s_1(i\zeta_1\tau,x)-1)\psi_2^*(i\tau,x)-(s_2(i\zeta_2\tau,x)-1)\psi_1^*(i\tau,x)\}\frac{d\tau}{\tau-i\zeta_1\widehat{\varkappa}_s}}\\
(1\leq s\leq\widehat{m}).
\end{array}\right.\label{eq2.61}
\end{equation}

\begin{conclusion}
Equalities \eqref{eq2.58}, \eqref{eq2.61} form the {\bf closed system} of linear singular equations relative to $\psi_1^*(it,x)$, $\psi_2^*(it,x)$ and $\{R_l(x)\}_1^m$, $\{\widehat{R}_s(x)\}_1^{\widehat{m}}$. $sc_1(\lambda)$, $sc_2(\lambda)$ and the number sets $\{\varkappa_l\}_1^m$, $\{\widehat\varkappa_s\}_1^{\widehat{m}}$ from $\mathbb{R}_+$ and $\{b_l\}_1^m$, $\{\widehat{b}_s\}_1^{\widehat{m}}$ from $\mathbb{C}$ are free (independent) variables of this system.
\end{conclusion}

\section{Inverse problem}\label{s3}

{\bf 3.1} Scheme of solution of an inverse problem is as follows. Let $\mathfrak{A}$ be a totality
\begin{equation}
\mathfrak{A}\stackrel{\rm def}{=}\{sc_1(\lambda);sc_2(\lambda);\{\varkappa_l,b_l\}_1^m;\{\widehat{\varkappa}_s,\widehat{b}_s\}_1^{\widehat{m}})\label{eq3.1}
\end{equation}
where $sc_1(\lambda)$, $sc_2(\lambda)$ are scattering coefficients \eqref{eq2.21}; $\{\varkappa_l\}_1^m$, $\{\widehat{\varkappa}_s\}_1^m$ are sets of positive numbers ($\varkappa_l$, $\varkappa_s>0$; $m$, $\widehat{m}\in\mathbb{Z}_+$); $\{b_l\}_1^m$, $\{\widehat{b}_s\}_1^{\widehat{m}}$ are complex numbers ($b_l$, $\widehat{b}_s\in\mathbb{C}$). From $sc_1(\lambda)$ and $sc_2(\lambda)$, construct $s_1(\lambda,x)$ \eqref{eq2.34} and $s_2(\lambda,x)$ \eqref{eq2.36} and afterwards consider the system of $2+m+\widehat{m}$ linear singular equations \eqref{eq2.58}, \eqref{eq2.61}. Knowing solution to this system $\psi_1^*(it,x)$, $\psi_2^*(it,x)$, $\{R_l\}_1^m$, $\{\widehat{R}_s(x)\}_1^{\widehat{m}}$, define the holomorphic in sector $\Omega_0^-$ function $\psi_0^*(\lambda,x)$ \eqref{eq2.57} and using \eqref{eq1.44}, we obtain
$$\int\limits_x^\infty q(t)dt=\lim\limits_{\lambda\rightarrow\infty}3i\lambda^2(\psi_0^*(\lambda,x)-1)\quad(\lambda\in\Omega_0^-;|\lambda|\gg1),$$
whence, due to \eqref{eq2.57}, we have
\begin{equation}
\begin{array}{ccc}
{\displaystyle\int\limits_x^\infty q(t)dt=3i\left\{-\zeta_1\sum\limits_1^m\varkappa_lR_l(x)+\zeta_2\sum\limits_1^m\widehat{\varkappa}_s\widehat{R}_s(x)+\frac1{2\pi i}\int\limits_0^\infty d\tau\left[(\zeta_0-\zeta_2)\cdot\int\limits_0^\tau dt\right.\right.}\\
{\displaystyle\times\left.\left.s_2(i\zeta_2t,x)\psi_1^*(it,x)+(\zeta_1-\zeta_0)\int\limits_0^\tau dts_1(i\zeta_1t,x)\psi_2^*(it,x)-\int\limits_0^\tau\Delta(t,x)dt\right]\right\};}
\end{array}\label{eq3.2}
\end{equation}
where $\Delta(t,x)\stackrel{\rm det}{=}\psi_2^*(it,x)-\psi_1^*(it,x)$. So, potential $q(x)$ is explicitly expressed \eqref{eq3.2} via the solution $\psi_1^*(it,x)$, $\psi_2^*(it,x)$, $\{R_l(x)\}_1^m$, $\{\widehat{R}_s(x)\}$ of system \eqref{eq2.58}, \eqref{eq2.61}.

In terms of
\begin{equation}
a(t,\lambda)\stackrel{\rm def}{=}\frac1{(t-\lambda)(t+\zeta_1\lambda)};\quad b(t,\lambda)\stackrel{\rm def}{=}\frac1{(t+\lambda)(t+\zeta_1\lambda)},\label{eq3.3}
\end{equation}
system \eqref{eq2.58}, considering Sokhotski formulas \cite{20,21}, becomes
\begin{equation}
\begin{array}{cccc}
{\displaystyle\psi_2^*(it,x)=1+\zeta_2\sum\limits_1^m\varkappa_lR_l(x)a(t,i\zeta_1\varkappa_l)-\zeta_0\sum\limits_1^{\widehat{m}}\widehat{\varkappa}_s\widehat{R}_s(x)a(t,i
\zeta_2
\widehat{\varkappa}_s)}\\
{\displaystyle+\frac12s_2(i\zeta_2t,x)\psi_1^*(it,x)+\frac{(\zeta_2-\zeta_1)t}{2\pi i}\int\limits_0^\infty\hspace{-4.4mm}/s_2^2(\lambda\zeta_2t,x)\psi_1^*(it,x)b(\tau,-t)}\\
{\displaystyle+\frac{(\zeta_1-\zeta_2)t}{2\pi i}\int\limits_0^\infty s_1(i\zeta_1\tau,x)\psi_2^*(i\tau,x)b(\tau,-\zeta_1t)d\tau-\frac1{2\pi i}\int\limits_0\infty\Delta(\tau,x)\frac{d\tau}{\tau-\zeta_2t};}\\
{\displaystyle\psi_1^*(it,x)=1+\zeta_0\sum\limits_1^m\varkappa_lR_l(x)a(t,i\varkappa_l)-\zeta_1\sum\limits_1^{\widehat{m}}\varkappa_s\widehat{R}_s(x)a(t,i\zeta_1
\widehat{\varkappa}_s)}
\end{array}\label{eq3.4}
\end{equation}
$$+\frac12s_1(i\zeta_1t,x)\psi_2^*(it,x)+\frac{(\zeta_2-\zeta_0)t}{2\pi i}\int\limits_0^\infty\hspace{-4.4mm}/s_1(i\zeta_1\tau,x)\psi_2^*(i\tau,x)b(\tau,-\zeta_2t)d\tau$$
$$+\frac{(\zeta_1-\zeta_2)t}{2\pi i}\int\limits_0^\infty s_2(i\zeta_2\tau,x)\psi_1^*(it,x)b(\tau,-\zeta_1t)dt-\frac1{2\pi i}\int\limits_0^\infty\Delta(\tau,x)\frac{d\tau}{\tau-\zeta_2t};$$
and, correspondingly, system \eqref{eq2.61} becomes
\begin{equation}
%\left\{
\begin{array}{ccccccccc}
{\displaystyle c_le^{i\varkappa_l(\zeta_2-\zeta_0)x}\varkappa_lR_l(x)=1+\zeta_2\sum\limits_1^m\varkappa_pR_p(x)a(\varkappa_p,\zeta_2\varkappa_l)}\\
{\displaystyle-\zeta_1\sum\limits_1^{\widehat{m}}\widehat{
\varkappa}_s\widehat{R}_s(x)a(\widehat{\varkappa}_s,\zeta_1\varkappa_l)}\\
{\displaystyle+\frac{(\zeta_2-\zeta_1)\varkappa_l}{2\pi}\int\limits_0^\infty s_2(i\zeta_2\tau,x)\psi_1^*(i\tau,x)b(\tau,i\zeta_1\varkappa_l)d\tau}\\
{\displaystyle-\frac1{2\pi i}\int\limits_0^\infty\Delta(\tau,x)\frac{d\tau}{\tau+i\zeta_2\varkappa_l}\quad(1\leq l\leq m, c_l\stackrel{\rm def}{=}b_l/\varkappa_l);}\\
{\displaystyle\widehat{c}_se^{-i\widehat{\varkappa}_s(\zeta_1-\zeta_0)x}\widehat{\varkappa}_s\widehat{R}_s(x)=1}\\
{\displaystyle+\zeta_2\sum\varkappa_lR_l(x)a(\varkappa_l,\zeta_2\widehat{
\varkappa}_s)+\zeta_1\sum\limits_1^{\widehat{m}}\widehat{\varkappa}_p\widehat{R}_p(x)a(\widehat{\varkappa}_p,\zeta_1\widehat{\varkappa}_s)}\\
{\displaystyle+\frac{(\zeta_0-\zeta_1)\widehat{\varkappa}_s}{2\pi}\int\limits_0^\infty s_2(i\zeta_2\tau,x)\psi_1^*(i\tau,x)b(\tau,-i\zeta_2\widehat{\varkappa}_s)d\tau}\\
{\displaystyle +\frac{(\zeta_0-\zeta_2)\widehat{\varkappa}_s}{2\pi}\int\limits_0^\infty s_1(i\zeta_1\tau,x)\psi_2^*(i\tau,x)b(\tau,-i\zeta_1\widehat{\varkappa}_s)d\tau}\\
{\displaystyle-\frac1{2\pi i}\int\limits_0^\infty\Delta(\tau,x)\frac{d\tau}{\tau-i\zeta_1\widehat{\varkappa}_s}\quad(1\leq s\leq\widehat{m},\widehat{c}_s=\widehat{b}_s/\widehat{\varkappa}_s).}
\end{array}\label{eq3.5}
\end{equation}

{\bf 3.2} First, consider the {\bf reflectionless potential} assuming that $sc_1(\lambda)=0$ and $sc_2(\lambda)=0$. Then system \eqref{eq3.4} becomes
$$\left\{
\begin{array}{cccc}
{\displaystyle\psi_2^*(it,x)=1+\zeta_2\sum\limits_1^m\varkappa_lR_l(x)a(t,i\zeta_1\varkappa_l)-\zeta_0\sum\limits_1^m\widehat{\varkappa}_s\widehat{R}_s(x)a(t,i\zeta_2
\widehat{\varkappa}_s)}\\
{\displaystyle-\frac1{2\pi i}\int\limits_0^\infty\Delta(\tau,x)\frac{d\tau}{\tau-\zeta_2t};}\\
{\displaystyle\psi_1^*(it,x)=1+\zeta_0\sum\limits_1^\infty\varkappa_lR_l(x)a(t,i\varkappa_0)-\zeta_1\sum\limits_1^{\widehat{m}}\widehat{\varkappa}_s\widehat{R}_s(x)a(t,i\zeta_1
\widehat{\varkappa}_s)}\\
{\displaystyle-\frac1{2\pi i}\int\limits_0^\infty\Delta(\tau,x)\frac{d\tau}{\tau-\zeta_2t}}.
\end{array}\right.$$
Subtracting the second equality from the first, we obtain
\begin{equation}
\Delta(t,x)=c(t,x)+M\Delta(t,x)\label{eq3.6}
\end{equation}
where
\begin{equation}
c(t,x)\stackrel{\rm def}{=}\sum\limits_1^m\varkappa_lR_l(x)A_l(t)+\sum\widehat{\varkappa}_s\widehat{R}_s(x)\widehat{A}_s(t);\label{eq3.7}
\end{equation}
\begin{equation}
A_l(t)\stackrel{\rm def}{=}\zeta_2a(t,i\zeta_1\varkappa_l)-\zeta_0a(t,i\varkappa_l);\quad\widehat{A}_s(t)\stackrel{\rm def}{=}\zeta_1a(t,i\zeta_1\widehat{\varkappa}_s)-\zeta_0a(t,i\zeta_2\widehat{\varkappa}_s);\label{eq3.8}
\end{equation}
and the Fredholm operator $M$ equals
\begin{equation}
(Mf)(t)\stackrel{\rm def}{=}\frac{(\zeta_2-\zeta_1)t}{2\pi i}\int\limits_0^\infty f(\tau)b(\tau,-\zeta_1t)d\tau\quad(f\in L^2(\mathbb{R}_+))\label{eq3.9}
\end{equation}
($b(t,\lambda)$ is from \eqref{eq3.3}). Equation \eqref{eq3.6} implies that $\Delta(t,x)=(1-M)^{-1}c(t,x)$ and, due to \eqref{eq3.7}, we have
\begin{equation}
\Delta(t,x)=\sum\limits_1^m\varkappa_lR_l(x)(1-M)^{-1}A_l(t)+\sum\limits_1^{\widehat{m}}\widehat{\varkappa}_s\widehat{R}_s(x)(1-M)^{-1}\widehat{A}_s(t).\label{eq3.10}
\end{equation}
Upon substituting this expression in system \eqref{eq3.5} in which $s_1(\lambda,x)=s_2(\lambda,x)=0$, we obtain
\begin{equation}
\left\{
\begin{array}{ccccccc}
{\displaystyle c_le^{i\varkappa_l(\zeta_2-\zeta_0)x}\varkappa_lR_l(x)=1+\sum\limits_1^m\varkappa_pR_p(x)\left[\zeta_2a(\varkappa_p,\zeta_2\varkappa_l)-\frac1{2\pi i}\int\limits_0^\infty(I-M)^{-1}\right.}\\
{\displaystyle \left.A_p(\tau)\cdot\frac{d\tau}{\tau+i\zeta_2\varkappa_l}\right]}\\
{\displaystyle+\sum\limits_1^{\widehat{m}}\widehat{\varkappa}_s\widehat{R}_s(x)\left[-\zeta_1a(\widehat{\varkappa}_s,\zeta_1\varkappa_l)-\frac1{2\pi i}\int\limits_0^\infty(I-M)^{-1}\widehat{A}_s(\tau)\frac{d\tau}{\tau+i\zeta_2\varkappa_l}\right]\quad(1\leq l\leq m);}\\
{\displaystyle\widehat{c}_se^{-i\widehat{\varkappa}_s(\zeta_1-\zeta_0)x}\widehat{\varkappa}_s\widehat{R}_s(x)=1+\sum\limits_1^m\varkappa_lR_l(x)\left[\zeta_2a(\varkappa_l,
\zeta_2\widehat{\varkappa}_s)-\frac1{2\pi i}\int\limits_0^\infty(I-M)^{-1}\right.}\\
{\displaystyle\left.A_l(\tau)\frac{d\tau}{\tau-i\zeta_1\widehat{\varkappa}_s}\right]}\\
{\displaystyle+\sum\limits_1^{\widehat{m}}\varkappa_p\widehat{R}_p(x)\left[\zeta_1a(\widehat{\varkappa}_p,\zeta_1\widehat{\varkappa}_s)-\frac1{2\pi i}\int\limits_0^\infty(I-M)^{-1}\widehat{A}_p(\tau)\frac{d\tau}{\tau-i\zeta_1\widehat{\varkappa}_s}\right]\quad(1\leq s\leq\widehat{m}).}
\end{array}\right.\label{eq3.11}
\end{equation}
Note that coefficients in \eqref{eq3.11} depend only on $\{\varkappa_l,b_l\}_1^m$ and $\{\widehat{\varkappa}_s,\widehat{b}_s\}_1^{\widehat{m}}$. Using \eqref{eq3.10}, rewrite formula \eqref{eq3.2} as
\begin{equation}
\begin{array}{ccc}
{\displaystyle\int\limits_x^\infty q(t)dt=3i\left\{\sum\limits_1^m\varkappa_lR_l(x)\left[-\zeta_1-\frac1{2\pi i}\int\limits_0^\infty\int\limits_0^\tau(I-M)^{-1}A_l(t)dtd\tau\right]\right.}\\
{\displaystyle\left.+\sum\limits_1^{\widehat{m}}\widehat{\varkappa}_s\widehat{R}_s\left[\zeta_2-\frac1{2\pi i}\int\limits_0^\infty\int\limits_0^\tau(I-M)^{-1}\widehat{A}_s(\tau)d\tau dt\right]\right\}}
\end{array}\label{eq3.12}
\end{equation}
where $\{A_l(t)\}_1^m$ and $\{\widehat{A}_s(t)\}_1^{\widehat{m}}$ are from \eqref{eq3.8}.

\begin{theorem}\label{t3.1}
If $sc_1(\lambda)=0$ and $sc_2(\lambda)=0$, then potential $q(x)$ is calculated from the solution $\{\varkappa_lR_l(x)\}_1^m$, $\{\widehat{\varkappa}_s\widehat{R}_s(x)\}$ to system \eqref{eq3.11} via formula \eqref{eq3.12} where $A_l(t)$, $\widehat{A}_s(t)$ and $M$ are given by \eqref{eq3.8} and \eqref{eq3.9} correspondingly.
\end{theorem}

{\bf Example.} Let $m=1$ ($\varkappa_1=\varkappa>0$) and $\widehat{m}=0$, then from the first equation in \eqref{eq3.11} we find that
$$\varkappa h(x)=\frac1{b(\varkappa)e^{i\varkappa(\zeta_2-\zeta_0)x}+cc(\varkappa)}\quad(R(x)=R_1(x))$$
where
$$b(\varkappa)=b_1\varkappa^{-1};\quad c(\varkappa)=\frac{\zeta_2}{(1-\zeta_2)\varkappa^2}-\frac1{2\pi i}\int\limits_0^\infty(I-M)^{-1}A_1(t)\frac{dt}{t+i\zeta_2\varkappa},$$
besides, $M$ and $A_1(t)$ are given by \eqref{eq3.9} and \eqref{eq3.2}. Taking \eqref{eq3.12} into account, we obtain
$$\int\limits_x^\infty q(t)=\varkappa R(x)a(\varkappa);\quad a(\varkappa)=-3i\zeta_1-\frac3{2\pi}\int\limits_0^\infty\int\limits_0^\tau(I-M)^{-1}A_1(t)dtd\tau,$$
and thus
$$q(x)=\frac{a(\varkappa)i(\zeta_2-\zeta_0)e^{i\varkappa(\zeta_2-\zeta_0)x}}{(b(x)e^{i\varkappa(\zeta_2-\zeta_0)x}+c(\varkappa))^2}$$
$$=\frac{i\varkappa(\zeta_2-\zeta_0)a(\varkappa)}
{\displaystyle(b(\varkappa)\cdot\exp\left(i\frac\varkappa2(\zeta_2-\zeta_0)x\right)+c(\varkappa)\exp\left(i\frac\varkappa2(\zeta_2-\zeta_0)x)\right)^2}.$$
This is analogous to the reflectionless potential (``soliton'') for the Schr\"{o}dinger operator.
\vspace{5mm}

{\bf 3.3} Suppose that $sc_2(\lambda)=0$ ($s_2(\lambda,x)=0$). Subtracting equations \eqref{eq3.4}, we obtain
$$\Delta(t,x)=c(t,x)-\frac12s_1(i\zeta_1t,x)\psi_2^*(it,x)+\frac{(\zeta_1-\zeta_2)t}{2\pi i}\int\limits_0^\infty s_1(i\zeta_2\tau,x)\psi_2^*(i\tau,x)b(\tau,-\zeta_1t)d\tau$$
$$+\frac{(\zeta_0-\zeta_2)t}{2\pi i}\int\limits_0^\infty\hspace{-4.4mm}/s_1(i\zeta_1\tau,x)\psi_1^*(i\tau,x)b(\tau,-\zeta_2t)d\tau+\frac{(\zeta_2-\zeta_1)t}{2\pi i}\int\limits_0^\infty\Delta(\tau,x)b(\tau,-\zeta_1t)d\tau$$
where $c(t,x)$ equals \eqref{eq3.7}. This equality and the first equation in \eqref{eq3.4} form the system of equations for determining $s_1(i\zeta_1t,x)\psi_2^*(it,x)$ and $\Delta(t,x)$. Rewrite this system in a matrix form:
\begin{equation}
T(t,x)\Phi(t,x)=h(t,x)+{\mathcal M}\Phi(t,x),\label{eq3.13}
\end{equation}
here $\Phi(t,x)$ and $h(t,x)$ are vector functions
\begin{equation}
\begin{array}{lll}
\Phi(t,x)\stackrel{\rm def}{=}\col[s_1(i\zeta_1t,x)\psi_2^*(it,x),\Delta(t,x)];\\
{\displaystyle h(t,x)\stackrel{\rm def}{=}\col\left[1+\zeta_2\sum\limits_1^m\varkappa_lR_l(x)a(t,i\zeta_1\varkappa_l)-\zeta_0\sum\limits_1^{\widehat{m}}\widehat{\varkappa}_s\widehat{R}_s(x)a(t,i\zeta_2\widehat{
\varkappa}_s,c(t,x))\right]}
\end{array}\label{eq3.14}
\end{equation}
($c(t,x)$ is given by \eqref{eq3.7}, \eqref{eq3.8}); matrix $T(t,x)$ is
\begin{equation}
T(t,x)\stackrel{\rm def}{=}\left[
\begin{array}{ccc}
s_1^{-1}(i\zeta_1t,x)&0\\
{\displaystyle\frac12}&1
\end{array}\right];\label{eq3.15}
\end{equation}
and operator ${\mathcal M}$ is given by
\begin{equation}
(\mathcal{M}f)(t)\stackrel{\rm def}{=}\frac1{2\pi i}\int\limits_0^mM(t,\tau)f(\tau)d\tau\quad(f\in L^2(\mathbb{R}_+,E^2)),\label{eq3.16}
\end{equation}
the kernel of which is matrix
\begin{equation}
M(t,\tau)\stackrel{\rm def}{=}\left[
\begin{array}{ccc}
(\zeta_1-\zeta_2)tb(\tau,-\zeta_1t)&-(\tau-\zeta_1t)^{-1}\\
(\zeta_1-\zeta_2)tb(\tau,\zeta_1t)+(\zeta_0-\zeta_2)tb(\tau,-\zeta_2t)&(\zeta_2-\zeta_1)tb(\tau,\zeta_1t)
\end{array}\right].\label{eq3.17}
\end{equation}
Since
$$\widetilde{T}(t,x)\stackrel{\rm def}{=}T^{-1}(t,x)=\left[
\begin{array}{ccc}
s_1(i\zeta_1t,x)&0\\
{\displaystyle-\frac12s_1(i\zeta_1t,x)}&1
\end{array}\right],$$
then, assuming
\begin{equation}
\widetilde{h}(t,x)\stackrel{\rm def}{=}\widetilde{T}(t,x)h(t,x);\quad\widetilde{M}(x,t,\tau)\stackrel{\rm def}{=}\widetilde{T}(t,x)M(t,\tau),\label{eq3.18}
\end{equation}
we find solution to equation \eqref{eq3.13},
\begin{equation}
\Phi(t,x)=(I-\widetilde{\mathcal{M}})^{-1}\widetilde{h}(t,x)\label{eq3.19}
\end{equation}
where $\widetilde{\mathcal M}$ is the integral operator \eqref{eq3.16} with kernel $\widetilde{M}(x,t,\tau)$ \eqref{eq3.18}.

Substitute the expression for $\Phi(t,x)$ \eqref{eq3.19} into equations of system \eqref{eq3.5}
\begin{equation}
\left\{
\begin{array}{cccc}
{\displaystyle c_le^{i\varkappa_l(\zeta_2-\zeta_1)x}\varkappa_lR_l(x)=1+\zeta_2\sum\limits_1^m\varkappa_p,\zeta_2\varkappa_l)-zeta_1\sum\limits_1^{
\widehat{m}}\widehat{\varkappa}_s\widehat{R}_s(x)a(\widehat{\varkappa}_s,\zeta_1\varkappa_l)}\\
{\displaystyle+\frac1{2\pi i}\int\limits_0^\infty\langle(I-\widetilde{M})^{-1}\widetilde{f}_1(\tau,x),f_l(\tau)\rangle_{E^2}d\tau\quad(1\leq l\leq m);}\\
{\displaystyle\widehat{c}_se^{-i\widehat{\varkappa}_s(\zeta_1-\zeta_0)x}\widehat{\varkappa}_s\widehat{R}_s(x)=1+\zeta_2\sum\limits_1^m\varkappa_lR_l(x)a(\varkappa_l,\zeta_2
\widehat{\varkappa}_s)+\zeta_1\sum\limits_1^{\widehat{m}}\widehat{\varkappa}_p\widehat{R}_p(x)a(\widehat{\varkappa}_p,\widehat{\zeta}_1\varkappa_s)}\\
{\displaystyle+\int\limits_0^\infty\langle(I-\widetilde{\mathcal{M}})^{-1}\widetilde{f}_1(\tau,x),\widehat{f}_s(\tau)\rangle_{E^2}d\tau\quad(1\leq s\leq\widehat{m}),}
\end{array}\right.\label{eq3.20}
\end{equation}
besides,

\begin{equation}
\begin{array}{lll}
f_l(t)\stackrel{\rm def}{=}\col[\overline{(\zeta_0-\zeta_2)i\varkappa_lb(t,i\zeta_2\varkappa_0)},-\overline{(t+i\zeta_2\varkappa_l)^{-1}}]\quad(1\leq l\leq m);\\
\widehat{f}_s(t)\stackrel{\rm def}{=}\col[\overline{(\zeta_1-\zeta_2)i\widehat{\varkappa}_sb(t,-i\zeta_1\widehat{\varkappa}_s)},-\overline{(t-i\zeta_1\varkappa_s)^{-1}}]\quad(1\leq s\leq\widehat{m}).
\end{array}\label{eq3.21}
\end{equation}
In this case, formula \eqref{eq3.2} becomes
\begin{equation}
\begin{array}{cccc}
{\displaystyle\int\limits_x^\infty q(t)dt=3i\left\{-\zeta_1\sum\limits_1^\infty\varkappa_lR_l(x)+\zeta_2\sum\limits_1^{\widehat{m}}\widehat{\varkappa}_s\widehat{R}_s(x)\right.}\\
{\displaystyle\left.+\frac1{2\pi i}\int\limits_0^\infty d\tau\int\limits_0^\tau\langle(I-\widetilde{\mathcal{M}})^{-1}\widetilde{h}(t,x),g\rangle dt\right\}}
\end{array}\label{eq3.22}
\end{equation}
where $g\stackrel{\rm def}{=}\col(\zeta_2-\zeta_0,-1)$.

\begin{theorem}\label{t3.2}
If $sc_2(\lambda)=0$ ($s_2(\lambda,x)=0$), then potential $q(x)$ is calculated from the solution $\{\varkappa_lR_l(x)\}_1^m$, $\{\widehat{\varkappa}_s\widehat{R}_s(x)\}$
to the system of equations \eqref{eq3.20} using formula \eqref{eq3.22} where $\widetilde{\mathcal{M}}$ and $\widetilde{h}(t,x)$ are given by \eqref{eq3.16} (with kernel $\widetilde{\mathcal M}(x,t,\tau)$ \eqref{eq3.18}) and \eqref{eq3.18}.
\end{theorem}

Inverse problem for arbitrary $sc_1(\lambda)\not=0$, $sc_2\not=0$ is solved analogously.

\renewcommand{\refname}{ \rm 
\end{Large}
\end{document}